\documentclass[11pt]{amsart}
\usepackage{amsmath, amssymb}
\usepackage{amsfonts}
\usepackage[arrow,matrix,curve,cmtip,ps]{xy}

\usepackage{amsthm}

\allowdisplaybreaks

\newtheorem{theorem}{Theorem}[section]
\newtheorem{lemma}[theorem]{Lemma}
\newtheorem{proposition}[theorem]{Proposition}
\newtheorem{corollary}[theorem]{Corollary}

\newtheorem*{theorem*}{Theorem}
\theoremstyle{remark}
\newtheorem{remark}[theorem]{Remark}
\newtheorem{definition}[theorem]{Definition}
\newtheorem{example}[theorem]{Example}

\numberwithin{equation}{section}

\newcommand{\Z}{\mathbb{Z}}
\newcommand{\N}{\mathbb{N}}
\newcommand{\C}{\mathbb{C}}
\newcommand{\T}{\mathbb{T}}

\newcommand{\F}{\mathcal{F}}
\newcommand{\G}{\mathcal{G}}
\newcommand{\La}{\mathcal{L}}

\newcommand{\reg}{\textnormal{reg}}
\newcommand{\infin}{\textnormal{inf}}

\newcommand{\algspan}{\operatorname{span}}

\begin{document}
\title[Leavitt path algebras]{Uniqueness Theorems and Ideal Structure for Leavitt Path Algebras}

\author{Mark Tomforde 
}

\address{Department of Mathematics \\ University of Houston \\ Houston, TX 77204-3008 \\USA}
\email{tomforde@math.uh.edu}


\date{\today}



\begin{abstract}

We prove Leavitt path algebra versions of the two uniqueness theorems of graph $C^*$-algebras.  We use these uniqueness theorems to analyze the ideal structure of Leavitt path algebras and give necessary and sufficient conditions for their simplicity.  We also use these results to give a proof of the fact that for any graph $E$ the Leavitt path algebra $L_\C(E)$ embeds as a dense $*$-subalgebra of the graph $C^*$-algebra $C^*(E)$.  This embedding has consequences for graph $C^*$-algebras, and we discuss how we obtain new information concerning the construction of $C^*(E)$.

\end{abstract}

\maketitle

\section{Introduction}

In the late 1950's Leavitt constructed examples of rings $R$ that do not have the ``Invariant Basis Number" property; i.e., ${}_RR^m \cong {}_RR^n$ as left $R$-module with $m \neq n$ \cite{Leav1, Leav2, Leav3}.  It can be shown that if $R$ is a unital ring, then ${}_RR^1 \cong {}_RR^n$ for $n > 1$ if and only if there exist elements $x_1, \ldots x_n, y_1, \ldots, y_n \in R$ such that $x_iy_j = \delta_{ij} 1_R$ for all $i,j$ and $\sum_{i=1}^n y_ix_i = 1_R$.  For a given field $K$, Leavitt considered the unital $K$-algebra $L(1, n)$ generated by elements $\{x_1, \ldots, x_n, y_1, \ldots, y_n \}$ satisfying these relations.

In 1977 Cuntz, independent of Leavitt's work, introduced a class of $C^*$-algebras generated by non-unitary isometries \cite{Cun}.  Specifically, if $n >1$ the Cuntz algebra $\mathcal{O}_n$ is the $C^*$-algebra generated by isometries $\{s_1, \ldots s_n \}$ satisfying $s_i^*s_j = \delta_{ij} I$ and $\sum_{i=1}^n s_is_i^* = I$.  The Cuntz algebras were the first examples of $C^*$-algebras that exhibited torsion in their $K$-theory \cite{Cun2} and have been shown to arise in many seemingly unrelated situations.  Consequently, the Cuntz algebras have been studied extensively, and their theory is now a standard part of the toolkit of $C^*$-algebraists.

A generalization of the Cuntz algebras was introduced in 1980  by Cuntz and Krieger \cite{CK}, where they considered $C^*$-algebras, now called Cuntz-Krieger algebras, associated to finite matrices with entries in $\{0, 1 \}$.  Approximately seventeen years later, in 1997, Kumjian, Pask, Raeburn and Renault \cite{KPR, KPRR} considered a generalization in which they associate a $C^*$-algebra to a (possibly infinite) directed graph.  These graph $C^*$-algebras have proven to be very important objects of study for $C^*$-algebraists.  

For a given graph $E$, the graph $C^*$-algebra $C^*(E)$ is the universal $C^*$-algebra generated by partial isometries satisfying relations determined by the graph $E$.  Remarkably, it has also been found that much of the structure of the $C^*$-algebra  $C^*(E)$ is reflected in the graph $E$.  This results in a beautiful and elegant theory in which one may translate $C^*$-algebraic properties into graph theoretic properties and vice versa.

Furthermore, graph $C^*$-algebras have proven useful for a number of reasons.  A short list of which include:  (1) the graph $C^*$-algebras include many known examples of $C^*$-algebras and thereby they help to organize many prior known theories under a common rubric, (2) the use of a directed graph as a tool allows one to translate many complicated $C^*$-algebra properties into graph theoretic properties that can be visualized and examined combinatorially, (3) many intricate $C^*$-algebra computations (e.g., for Ext and $K$-theory) can be reformulated as graph computations that are typically easier to deal with, (4) one can produce examples and (counterexamples) of $C^*$-algebras with given properties simply by producing directed graphs with the corresponding properties.

Initially, graph $C^*$-algebras were considered only for row-finite graphs (i.e., graphs in which each vertex emits at most a finite number of edges), and in fact it was unclear how to extend the definition of $C^*(E)$ to non-row-finite graphs.  Later, in 2000, it was determined how to appropriately define $C^*(E)$ for arbitrary graphs \cite{FLR}.  The $C^*$-algebras of non-row-finite graphs were found to include examples of $C^*$-algebras not included in the row-finite case, and it was also discovered that non-row-finite graph $C^*$-algebras exhibited behavior that was different from the row-finite case.  Consequently, descriptions of the $C^*$-algebraic properties of $C^*(E)$ in terms of the properties of $E$ had to be modified significantly from the row-finite case, and new techniques had to be developed to extend results from the row-finite case to the general case.  Nonetheless, this program has been fairly successful, and today there exists an extensive theory of graph $C^*$-algebras which has found applications to many other areas of operator algebra.

Inspired by the success of graph $C^*$-algebras, G.~Abrams and G.~Aranda-Pino sought to create algebraic analogues of the graph $C^*$-algebras that generalize the relationship between the Leavitt algebra $L(1,n)$ and the Cuntz algebra $\mathcal{O}_n$. In 2005, Abrams and Aranda-Pino defined the \emph{Leavitt path algebra} $L_K(E)$ to be the universal $K$-algebra generated by elements satisfying relations similar to those of the generators for the graph $C^*$-algebra of $E$ \cite{AbrPino}.  As with graph $C^*$-algebras, this initial definition was given only in the case that $E$ is a row-finite graph.  Abrams and Aranda-Pino established basic properties of these Leavitt path algebras, and soon after there was a flurry of activity as many other authors \cite{AbrPino2, AbrPino3, APM, APM2, AMP2, AraPar, PinParMol} investigated their structure and found applications to various topics in algebra.  As with graph $C^*$-algebras, it was found that much of the structure of the algebra $L_K(E)$ is reflected in the graph $E$.  Furthermore, and to many researchers' surprise, it was often the case that when certain graph properties of $E$ correspond to $C^*$-algebraic properties of $C^*(E)$, the same graph properties of $E$ correspond to the analogous algebraic properties of $L_K(E)$.  A few examples of this are the following: the conditions on $E$ that correspond to simplicity of $C^*(E)$ are the same as the conditions on $E$ that correspond to simplicity of $L_K(E)$; the conditions on $E$ that correspond to $C^*(E)$ being simple and purely infinite are the same as the conditions on $E$ that correspond to $L_K(E)$ being simple and purely infinite; and the gauge-invariant ideals of $C^*(E)$ correspond to saturated hereditary subsets of vertices of $E$ in much the same way that the graded ideals of $L_K(E)$ correspond to saturated hereditary subsets of $E$.  

What is even more astonishing is that neither the graph $C^*$-algebra results nor the Leavitt path algebra results are obviously logical consequences of the other.  Often different methods are used in the proofs, and moreover, neither result can easily be seen to imply the other.  This has left many researchers wondering exactly what the relationship is between the Leavitt path algebras and the graph $C^*$-algebras. 

We also mention that Abrams and Aranda-Pino have very recently (September, 2006) extended the definition of the Leavitt path algebra to include non-row-finite graphs \cite[Definition~1]{AbrPino3}.  As with graph $C^*$-algebras, this has resulted in new examples of Leavitt path algebras not included in the row-finite case, and many of these new algebras exhibit behavior different from the Leavitt path algebras of row-finite graphs.  At the time this paper was written, \cite{AbrPino3} is the only other paper which has proven results about Leavitt path algebras of non-row-finite graphs.

The purpose of this paper is to establish analogues for Leavitt path algebras of the two main uniqueness theorems of graph $C^*$-algebras, to use these theorems to deduce facts about the ideal structure of Leavitt path algebras, and to examine the relationship between Leavitt path algebras and graph $C^*$-algebras.   In doing this we hope to establish techniques that will be useful for future authors as they study Leavitt path algebras.  We will make no assumption of row-finiteness of our graphs, and will prove all results in the general case.  As described earlier, this will involve a more in-depth analysis and will require techniques different from those used in the row-finite case.

There are two main uniqueness theorems for graph $C^*$-algebras: The Gauge-Invariant Uniqueness Theorem and the Cuntz-Krieger Uniqueness Theorem.  Each of these theorems gives conditions under which homomorphisms on $C^*(E)$ are injective.  They are fundamental results in the subject and are used extensively in identifying the isomorphism class of a particular graph $C^*$-algebra as well as deducing many general results about the structure of all graph $C^*$-algebras.  Consequently, from the graph $C^*$-algebraists perspective, it is rather surprising that the algebraists who have worked on Leavitt path algebras have not explicitly stated and proven versions of these theorems.  Often in papers on Leavitt path algebras the authors have used ad hoc methods to prove that a particular homomorphism is injective, and in deeper results a version of the uniqueness theorem is often hidden in their proofs.  (E.g., \cite[Corollary~3.3]{AbrPino} can be interpreted as a weak version of the Cuntz-Krieger Uniqueness Theorem for Leavitt algebras of row-finite graphs; the proof of \cite[Theorem~5.3]{AMP2} involves showing that every graded ideal is generated by idempotents, which is essentially an analogue of the Gauge-Invariant Uniqueness Theorem in the row-finite case.)  In the author's opinion, it would be useful for the study of Leavitt path algebras if these uniqueness theorems were explicitly stated and proven, so that they can be used in a more forthright manner and be more accessible to other researchers.

After some definitions and preliminary results in \S\ref{def-sec} and \S \ref{basic-results-sec}, we proceed in \S\ref{GUT-sec} to prove a Graded Uniqueness Theorem for Leavitt path algebras.  This is an analogue of the Gauge-Invariant Uniqueness Theorem for graph $C^*$-algebras, but because the grading on $L_K(E)$ plays the role of the gauge action on $C^*(E)$, we change the name.  

In \S\ref{graded-ideals-sec} we use the Graded Uniqueness Theorem to analyze the ideal structure of a Leavitt path algebra, and to give a complete description of the lattice of graded ideals of $L_K(E)$ in terms of $E$.  This result generalizes \cite[Theorem~5.3]{AMP2}, and (as with graph $C^*$-algebras) we will need to consider new phenomena that arise in the non-row-finite case.  Whereas graded ideals of Leavitt path algebras of row-finite graphs correspond to saturated hereditary subsets of vertices of $E$, the graded ideals of a general Leavitt path algebra correspond to \emph{admissible pairs} $(H,S)$, where $H$ is a saturated hereditary set of vertices and $S$ is a special subset of vertices that emit infinitely many edges.  We also use the Graded Uniqueness Theorem to identify the quotient of a Leavitt path algebra by a graded ideal as a Leavitt path algebra, and to show that any graded ideal of a Leavitt path algebra is Morita equivalent to a Leavitt path algebra.

In \S\ref{CK-uniqueness-sec} we prove a version of the Cuntz-Krieger Uniqueness Theorem for Leavitt path algebras.  We use this theorem to show that, in analogy with the graph $C^*$-algebra theory, a graph $E$ satisfies Condition~(K) if and only if all ideals of $L_K(E)$ are graded.  We then use the Graded Uniqueness Theorem and the Cuntz-Krieger Uniqueness Theorem together to give necessary and sufficient conditions on $E$ for $L_K(E)$ to be simple.

Finally, in \S\ref{Leavitt-and-graph-algebras-sec} we look at the relationship between Leavitt path algebras and graph $C^*$-algebras.  It is easy to see that if $E$ is a graph, then the universal property of $L_K(E)$ implies there is a homomorphism $\phi : L_\C(E) \to C^*(E)$ taking generators to generators, and the image of $\phi$ is a dense $*$-subalgebra of $C^*(E)$.  It has been asserted by many authors that for row-finite graphs this homomorphism is injective, but a proof has never been written down in the literature.  We discuss some of the subtleties of this statement and argue that it is not immediately obvious.  (It is, of course, obvious when $L_\C(E)$ is simple --- so $L(1,n)$ embeds into $\mathcal{O}_n$, for example --- but it is not clear when $L_\C(E)$ is not simple.)  Nonetheless, we are able to show that the injectivity of $\phi$ follows from the Graded Uniqueness Theorem we prove in \S\ref{GUT-sec}.  Thus the result is true, but seemingly nontrivial.  We also show that the injectivity of $\phi$ has consequences for the construction of $C^*(E)$.  Typically $C^*(E)$ is constructed by taking a free algebra $k_E$ with generators corresponding to edges and vertices of $E$, forming a quotient $k_E /J$ to ensure the necessary relations hold, and then defining a semi-norm on $k_E/J$.  One then defines $K$ to be all elements in $k_E / J$ for which this semi-norm equals zero, and $C^*(E)$ is equal to the completion $\overline{(k_E/J)/K}$.  We show in Corollary~\ref{graph-C*-construct-cor} and Corollary~\ref{norm-LC-in-C*} that this semi-norm is actually a norm, so that $K= \{ 0 \}$ and $C^*(E) = \overline{k_E / J}$.  In addition, we show that by viewing $L_\C(E)$ as a $*$-subalgebra of $C^*(E)$, we obtain a correspondence $I \mapsto \overline{I}$ between graded ideals of $L_\C(E)$ and gauge-invariant ideals of $C^*(E)$.

After this paper was written, it was brought to the author's attention that Iain Raeburn had proven special cases of the uniqueness theorems for Leavitt path algebras \cite[\S1.3]{Rae2}.  In his work, however, there are hypotheses on the underlying field $K$ (viz.~that $K$ is a $*$-field, and that the $*$-operation is ``positive definite") as well as a standing hypothesis that all graphs are row-finite.  In the uniqueness theorems of this paper, Theorem~\ref{GUT} and Theorem~\ref{CK-Uniqueness}, we are able to circumvent the need for any of these hypotheses.

$ $

\noindent \textbf{Notation and Conventions:}  Since our audience is both algebraists and $C^*$-algebraists, we do our best to make this paper accessible to both groups.  Because we are working with rings and $K$-algebras throughout this paper, we discuss preliminaries and establish basic facts in \S\ref{basic-results-sec}.  While much of this may seem pedantic to our algebraist readers, it will be of more interest to $C^*$-algebraists since we focus on those aspects of rings and algebras that are different from what occurs in the $C^*$-algebra setting.  Throughout this paper we will also make attempts to compare our results for Leavitt path algebras to the graph $C^*$-algebra theory.  

In addition,  we will use the following conventions:   all of our rings are not necessarily commutative and do not necessarily have identity; all ideals are two-sided; our graphs are not necessarily row-finite; the symbol $K$ denotes a field; $M_n(K)$ denotes the $n \times n$ matrices with entries in $K$; and $M_\infty(K)$ denotes the infinite $\N \times \N$ matrices with all but a finite number of entries equal to zero  --- thus if we include $M_n(K) \hookrightarrow M_{n+1}(K)$ by mapping to the upper left corner, we have $M_\infty (K) = \bigcup_{n=1}^\infty M_n (K)$.

\section{The definition of the Leavitt path algebra} \label{def-sec}

In this paper when we refer to a graph, we shall always mean a directed graph $E := (E^0, E^1, r, s)$ consisting of a countable set of vertices $E^0$, a countable set of edges $E^1$, and maps $r: E^1 \to E^0$ and $s:E^1 \to E^0$ identifying the range and source of each edge.

\begin{definition}
Let $E := (E^0, E^1, r, s)$ be a graph.  We say that a vertex $v \in E^0$ is a \emph{sink} if $s^{-1}(v) = \emptyset$, and we say that a vertex $v \in E^0$ is an \emph{infinite emitter} if $|s^{-1}(v)| = \infty$.  A \emph{singular vertex} is a vertex that is either a sink or an infinite emitter, and we denote the set of singular vertices by $E^0_\textnormal{sing}$.  We also let $E^0_\textnormal{reg} := E^0 \setminus E^0_\textnormal{sing}$, and refer to the elements of $E^0_\textnormal{reg}$ as \emph{regular vertices}; i.e., a vertex $v \in E^0$ is a regular vertex if and only if $0 < |s^{-1}(v)| < \infty$.
\end{definition}

\begin{definition}
If $E$ is a graph, a \emph{path} is a sequence $\alpha := e_1 e_2 \ldots e_n$ of edges with $r(e_i) = s(e_{i+1})$ for $1 \leq i \leq n-1$.  We say the path $\alpha$ has \emph{length} $| \alpha| :=n$, and we let $E^n$ denote the set of paths of length $n$.  We consider the vertices in $E^0$ to be paths of length zero.  We also let $E^* := \bigcup_{n=0}^\infty E^n$ denote the paths of finite length, and we extend the maps $r$ and $s$ to $E^*$ as follows: For $\alpha := e_1 e_2 \ldots e_n \in E^n$, we set $r(\alpha) = r(e_n)$ and $s(\alpha) = s(e_1)$. 
\end{definition}

\begin{definition}
We let $(E^1)^*$ denote the set of formal symbols $\{ e^* : e \in E^1 \}$, and for $\alpha = e_1 \ldots e_n \in E^n$ we define $\alpha^* := e_n^* e_{n-1}^* \ldots e_1^*$.  We also define $v^* = v$ for all $v \in E^0$.  We call the elements of $E^1$ \emph{real edges} and the elements of $(E^1)^*$ ghost edges.
\end{definition}

\begin{definition} \label{Leavitt-def}
Let $E$ be a directed graph, and let $K$ be a field.  The \emph{Leavitt path algebra of $E$ with coefficients in $K$}, denoted $L_K(E)$,  is the universal $K$-algebra generated by a set $\{v : v \in E^0 \}$ of pairwise orthogonal idempotents, together with a set $\{e, e^* : e \in E^1\}$ of elements satisfying
\begin{enumerate}
\item $s(e)e = er(e) =e$ for all $e \in E^1$
\item $r(e)e^* = e^* s(e) = e^*$ for all $e \in E^1$
\item $e^*f = \delta_{e,f} \, r(e)$ for all $e, f \in E^1$
\item $v = \displaystyle \sum_{\{e \in E^1 : s(e) = v \}} ee^*$ whenever $v \in E^0_\reg$.
\end{enumerate}
\end{definition}

\begin{remark}
When we say that $L_K(E)$ is universal, we mean that it is universal for the relations listed in the definition.  In other words, if $A$ is a $K$-algebra containing a set of pairwise orthogonal idempotents $\{a_v : v \in E^0\}$ and a set of elements $\{ b_e, b_{e^*} : e \in E^1 \}$ satisfying the relations listed above, then there exists an algebra homomorphism $\phi : L_K(E) \to A$ with $\phi(v) = a_v$ for all $v \in E^0$ as well as $\phi(e) = b_e$ and $\phi(e^*) = b_{e^*}$ for all $e \in E^1$.
\end{remark}

\begin{remark}[The Construction of $L_K(E)$]  To show that $L_K(E)$ exists, begin by constructing the free algebra $K[E^0 \cup E^1 \cup (E^1)^*]$ subject to the relations
\begin{itemize}
\item[(i)] $vw = \delta_{vw} v \text{ for every $v,w \in E^0$}$
\item[(ii)] $e = e r(e) = s(e) e \text{ for every $e \in E^1$}$
\item[(iii)] $e^* = e^* s(e) = r(e) e \text{ for every $e \in E^1$.}$
\end{itemize}
(As pointed out in \cite[Definitions~1.2 and 1.3]{AbrPino}, this is the path algebra of the extended graph $\widehat{E}$ formed from $E$ by adding a ``ghost edge" $e^*$ for every $e \in E^1$ with $s(e^*) = r(e)$ and $r(e^*) = s(e)$.)  We then let $I$ be the ideal in this algebra generated by the elements $\{ e^*f - \delta_{ef} r(e) : e, f \in E^1 \} \cup \{ v - \sum_{s(e)=v} ee^* : v \in E^0_\textnormal{reg} \}$.  We define $L_K(E)$ to be the quotient of this algebra by the ideal $I$.  We see that in $L_K(E)$ the following relations hold:
\begin{itemize}
\item[(CK1)] $e^*f = \delta_{ef} r(e) \text{ for $e, f \in E^1$}$
\item[(CK2)] $v = \sum_{s(e)=v} ee^* \text{ for $v \in E^0_\textnormal{reg}$}.$ 
\end{itemize}
The relations (CK1) and (CK2) are often called the ``Cuntz-Krieger relations" after the relations introduced by Cuntz and Krieger for their eponymous $C^*$-algebras.  We can see that $L_K(E)$ as defined here satisfies the relations in Definition~\ref{Leavitt-def}, and has the appropriate universal property.  Furthermore, it can be shown that the elements $\{ v : v \in E^0 \} \cup \{ e, e^* : e \in E^1 \}$ are all nonzero in $L_K(E)$.  Finally, we mention that the universal property implies that $L_K(E)$ is unique (up to isomorphism).
\end{remark}

\section{Basic results for Leavitt path algebra} \label{basic-results-sec}

In this section we establish some basic results for Leavitt path algebras that will later be useful for us.  We shall try to compare and contrast these results to the corresponding facts for $C^*$-algebras.

\subsection{Involution and self-adjoint ideals}

We see that elements of $L_K(E)$ will be linear combinations of words in $\{ v : v \in E^0 \} \cup \{ e, e^* : e \in E^1 \}$ with coefficients from the field $K$.  From the relations in Definition~\ref{Leavitt-def} can see the following.

\begin{lemma} \label{CK-multiplication-lem}
If $E$ is a graph and $L_K(E)$ is the associated Leavitt path algebra, then for any $\alpha, \beta, \gamma, \delta \in E^*$ we have
$$(\alpha \beta^*) (\gamma \delta^*) = \begin{cases} \alpha \gamma' \delta^* & \text{ if $\gamma = \beta \gamma'$} \\ \alpha  \delta^* & \text{ if $\beta = \gamma$} \\ \alpha \beta'^* \delta^* & \text{ if $\beta = \gamma \beta'$} \\ 0 & \text{ otherwise.} \end{cases}$$
\end{lemma}

\noindent This result shows us that any word in $\{ v : v \in E^0 \} \cup \{ e, e^* : e \in E^1 \}$ may be written in the form $\alpha \beta^*$, where $\alpha, \beta \in E^*$ and $r(\alpha) = r(\beta)$.

\begin{corollary} \label{L(E)-span-cor}
If $E$ is a graph and $L_K(E)$ is the associated Leavitt path algebra, then $$L_K(E) = \algspan \{ \alpha \beta^* : \alpha, \beta \in E^*\text{ and } r(\alpha) = r(\beta) \}.$$  Consequently, any element $x \in L_K(E)$ can be written in the form $x = \sum_{k=1}^n \lambda_k \alpha_k \beta_k^*$, where $\lambda_k \in K$ and $\alpha_k, \beta_k \in E^*$ with $r(\alpha_k) = r(\beta_k)$ for $1 \leq k \leq n$.
\end{corollary}

\begin{definition}
We say that $x \in L_K(E)$ is a \emph{polynomial in all real edges} if $x = \sum_{k=1}^n \lambda_k \alpha_k$ for $\lambda_k \in K$ and $\alpha_k \in E^*$.  We say that $x \in L_K(E)$ is a \emph{polynomial in all ghost edges} if $x = \sum_{k=1}^n \lambda_k \beta^*_k$ for $\beta_k \in K$ and $\alpha_k \in E^*$.
\end{definition}

\begin{remark}
If $E$ is a graph and $L_K(E)$ is the associated Leavitt path algebra, we may define a linear involution $x \mapsto \overline{x}$ on $L_K(E)$ as follows:  If $x = \sum_{k=1}^n \lambda_k \alpha_k \beta_k^*$, then $\overline{x} = \sum_{k=1}^n \lambda_k \beta_k \alpha_k^*$.  Note that this operation is linear, involutive ($\overline{\overline{x}} = x$), and anti-multiplicative ($\overline{xy} = \overline{y} \, \overline{x}$).

When $K = \C$ we may also define a conjugate-linear involution $*$ on $L_\C(E)$ by setting $x^* := \sum_{k=1}^n \overline{\lambda_k} \beta_k^* \alpha_k$, where $\overline{\lambda_k}$ denotes the complex conjugate of $\lambda_k$.  This $*$-operation is conjugate-linear, involutive, and anti-multiplicative.
\end{remark}

\begin{definition}
If $L_K(E)$ is a Leavitt path algebra, an ideal $I$ of $L_K(E)$ is \emph{self-adjoint} if $\overline{I} = I$.
\end{definition}

\begin{remark} \label{non-selfadjoint-ideal-rem}
As described in \cite[Remark~3.5]{AbrPino}, unlike ideals in $C^*$-algebras, the ideals in Leavitt path algebras need not be self-adjoint.  If $E$ is the graph

$ $

$$
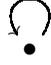
\xymatrix{
\bullet \ar@(ur,ul) \\
}
$$

\noindent consisting of a single vertex and a single edge, then $L_K(E) \cong K[x,x^{-1}]$.  If we let $p := 1 + x + x^3$, and $I:= \langle p \rangle$, then we can show $I$ is not self adjoint as follows:  If $\overline{p} = 1 + x^{-1} + x^{-3} \in I$, then $1 + x ^2+ x^3 = x^3 (1 + x^{-1} + x^{-3}) \in I$.  But then, since $K[x,x^{-1}]$ is commutative, we must have $q (1+x+x^3) = 1 + x^2 + x^3$ for some $q \in K[x,x^{-1}]$.  A degree argument then shows that $\deg q =0$, and $q = a_0 \in K$, which is absurd.  Hence $\overline{p} \notin I$, and $I$ is not self-adjoint.
\end{remark}

\subsection{Rings with local units}

The Leavitt path algebras that we look at will not necessarily have a unit.  However, every Leavitt path algebra does have a set of local units.  We mention that every $C^*$-algebra has an approximate unit, which can often play the role of a unit when a unit does not exist.  A set of local units plays a similar role in the ring setting.

\begin{definition}
A \emph{set of local units} for a ring $R$ is a set $E \subseteq R$ of commuting idempotents with the property that for any $x \in R$ there exists $t \in E$ such that $tx=xt=x$.
\end{definition}

By induction we obtain the following fact.

\begin{proposition}
If $R$ is a ring with a set of local units $E$, then for any finite number of elements $x_1, \ldots, x_n \in R$, there exists $t \in E$ such that $tx_i = x_it = x_i$ for all $1 \leq i \leq n$.
\end{proposition}

\begin{definition}[\cite{Ful}]
We say that a ring $R$ has \emph{enough idempotents} if there exists a collection of mutually orthogonal idempotents $\{ e_\alpha \}_{\alpha \in \Lambda}$ such that $R = \bigoplus e_\alpha R = \bigoplus R e_\alpha$.
\end{definition}

Note that if we let $S = \{ e_\alpha \}_{\alpha \in \Lambda}$ be the mutually orthogonal idempotents of the above definition, then $E := \{ \sum_{k=1}^n e_k : e_1, \ldots, e_n \in S\}$ is a set of local units for $R$.  Thus rings with enough idempotents are rings with local units.

\begin{remark}
If $E$ is a graph and $L_K(E)$ is the associated Leavitt path algebra, then $$L_K(E) = \bigoplus_{v \in E^0} v L_K(E) = \bigoplus_{v \in E^0} L_K(E) v$$ so $L_K(E)$ is a ring with enough idempotents.  Furthermore, if we list the vertices of $E$ as $E^0 = \{v_1, v_2, \ldots \}$, let $$\Lambda := \begin{cases} \{ 1, 2,  \ldots , |E^0| \} & \text{ if $E^0$ is finite} \\ \{1, 2, \ldots \} & \text{ if $E^0$ is infinite} \end{cases}$$
and set $t_n := \sum_{k=1}^n v_k$, then $\{ t_n \}_{n \in \Lambda}$ is a set of local units for $L_K(E)$.
\end{remark}

In general, if $A$ is a $K$-algebra, then a ring ideal of $A$ is not necessarily an algebra ideal of $A$.  However, when $A$ has a set of local units, the ring ideals and algebra ideals coincide.

\begin{lemma} \label{ring-ideals-alg-ideals-lem}
If $A$ is an algebra that is also a ring with a set of local units, then $I$ is a ring ideal of $R$ if and only if $I$ is an algebra ideal of $A$. 
\end{lemma}

\begin{proof}
If $I$ is an algebra ideal, then $I$ is trivially a ring ideal.  If $I$ is a ring ideal, then to show $I$ is an algebra ideal it suffices to show that $I$ is closed under scalar multiplication by elements of $K$.  Let $x \in I$ and $k \in K$.  Choose an idempotent $t \in R$ such that $tx = x$.  Since $I$ is a ring ideal, $kx = k(tx) = (kt) x \in I$.
\end{proof}

The above lemma will be useful when we discuss Morita equivalence of Leavitt path algebras.  Since the ring ideals and algebra ideals coincide, we will not have to worry about which we are discussing when we use the term ideal.

\subsection{$\Z$-graded rings}

All Leavitt path algebras have a natural $\Z$-grading.  As we shall see in \S \ref{GUT-sec}, this grading plays a role analogous to that of the gauge action for graph $C^*$-algebras.

\begin{definition}
If $R$ is a ring, we say $R$ is $\Z$-graded if there is a a collection of additive subgroups $\{ R_n \}_{n \in \Z}$ of $R$ with the following two properties.
\begin{enumerate}
\item $R = \bigoplus_{n \in \Z} R_n$.
\item $R_m R_n \subseteq R_{m+n}$ for all $m,n \in \Z$.
\end{enumerate}
The subgroup $R_n$ is called the \emph{homogeneous component of $R$ of degree $n$}.
\end{definition}

If $E$ is a graph, then we may define a $\Z$-grading on the associated Leavitt path algebra $L_K(E)$ by setting $$L_K(E)_n := \left\{ \sum_{k=1}^l \lambda_k \alpha_k \beta_k^* : \alpha_k, \beta_k \in E^* \text{ and } |\alpha_k | - | \beta_k| = n \text{ for $1 \leq k \leq l$} \right\}.$$
The fact that this is a grading follows from Lemma~\ref{CK-multiplication-lem} and Corollary~\ref{L(E)-span-cor}.  Note that, in fact, each $L_K(E)_n$ is closed under scalar multiplication by elements of $K$.  Hence $L_K(E)$ is actually a graded algebra.  However, in this paper we will be primarily concerned with the fact that $L_K(E)$ is a graded ring.

\begin{definition}
If $R$ is a graded ring, then an ideal $I$ of $R$ is a \emph{$\Z$-graded ideal} if $I = \bigoplus_{n \in \Z} (I \cap R_n)$.  If $\phi : R \to S$ is a ring homomorphism between $\Z$-graded rings, then $\phi$ is a \emph{graded ring homomorphism} if $\phi(R_n) \subseteq S_n$ for all $n \in \Z$.
\end{definition}

Note that the kernel of a $\Z$-graded homomorphism is a $\Z$-graded ideal.  Also, if $I$ is a $\Z$-graded ideal in a $\Z$-graded ring $R$, then the quotient $R / I$ admits a natural $\Z$-grading and the quotient map $R \to R /I$ is a $\Z$-graded homomorphism.  In this paper we will be concerned only with $\Z$-gradings, and hence we will often omit the prefix $\Z$ and simply refer to rings, ideals, homomorphisms, etc. ~as \emph{graded}.

\begin{definition}
A ring $R$ is \emph{idempotent} if $R^2 = R$; that is, if every $x \in R$ can be written as $x = \sum_{k=1}^n a_k b_k$ for $a_1, \ldots a_n, b_1, \ldots, b_n \in R$.
\end{definition}

\begin{remark}
We see that if $R$ is a ring with a set of local units, then $R$ is idempotent:  If $x \in R$, then there exists an idempotent $t \in R$ with $x = tx$.  Consequently, the Leavitt path algebra $L_K(E)$ is an idempotent ring.
\end{remark}

\subsection{Morita equivalence}

Throughout this paper we will need to discuss Morita equivalence for rings that do not necessarily have an identity element.  We establish the necessary definitions and results here.

\begin{definition}
If $R$ is a ring, we say that a left $R$-module is \emph{unital} if $RM = M$.  We also say that $M$ is \emph{non-degenerate} if for all $m \in M$ we have that $Rm = 0$ implies that $m = 0$.  We let $R-Mod$ denote the full subcategory of the category of all $R$-modules whose objects are unital non-degenerate $R$-modules.  (Note that if $R$ is unital, $R$-Mod is the usual category of $R$-modules.)  When $R$ and $S$ are rings, and ${}_RM_S$ is a bimodule, we say $M$ is \emph{unital} if $RM = M$ and $MS=M$.
\end{definition}

\begin{definition}
Let $R$ and $S$ be idempotent rings.  A \emph{(surjective) Morita context} $(R,S,M,N,\psi, \phi)$ between $R$ and $S$ consists of unital bimodules ${}_RM_S$ and ${}_SN_R$, a surjective $R$-module homomorphism $\psi : M \otimes_S N \to R$, and a surjective $S$-module homomorphism $\phi : N \otimes_R M \to S$ satisfying
$$\phi (n \otimes m) n' = n \psi (m \otimes n') \qquad \text{ and } \qquad m' \phi(n \otimes m) = \psi (m' \otimes n)m $$ for every $m, m' \in M$ and $n,n' \in N$.  We say that $R$ and $S$ are \emph{Morita equivalent} in the case that there exists a Morita context.
\end{definition}

It is proven in \cite[Proposition~2.5]{GS} and \cite[Proposition~2.7]{GS} that $R-Mod$ and $S-Mod$ are equivalent categories if and only if there exists a Morita context $(R,S,M,N,\psi, \phi)$.  In addition, the following result is obtained in \cite{GS}. 

\begin{proposition}  \cite[Proposition~3.5]{GS} \label{ideal-corresp-Mor-eq}
Let $R$ and $S$ be Morita equivalent idempotent rings, and let $(R,S,M,N,\psi, \phi)$ be a Morita context.  If $$\mathcal{L}_R := \{ I \subseteq R : \text{ $I$ is an ideal and $RIR=I$} \}$$ and $$\mathcal{L}_S := \{ I \subseteq S : \text{ $I$ is an ideal and $SIS=I$} \},$$ then there is a lattice isomorphism from $\mathcal{L}_R$ onto $\mathcal{L}_S$ given by $I \mapsto \phi (NI, M)$ with inverse given by $I \mapsto \psi (MI, N)$.
\end{proposition}

\begin{remark} \label{ideal-corresp-rings-loc-units}
Note that when $R$ is a ring with a set of local units, $\mathcal{L}_R$ is the lattice of ideals of $R$.  Similarly when $S$ is a ring with a set of local units.  Therefore if each of $R$ and $S$ is a ring with sets of local units, and if $R$ and $S$ are Morita equivalent, then the lattice of ideals of $R$ is isomorphic to the lattice of ideals of $S$.
\end{remark}

\begin{remark} \label{ideals-not-tran-remark}
Recall that, unlike $C^*$-algebras, the property of being a ring ideal is not transitive; i.e., if $R$ is a ring, $I$ is an ideal of $R$, and $J$ is an ideal of $I$, then it is not necessarily true that $J$ is an ideal of $R$.  (To see that this is the case, let $K$ be a field and let $R := K[x]$.  If we let $I := \{ p(x) \in R : a_0 = a_1 = 0 \}$ and let $J := \{ p(x) \in R : a_0=a_1=a_3=0 \}$, then $I \triangleleft R$, and $J \triangleleft I$.  However, $J$ is not an ideal of $R$ since $x^2 \in J$ and $x \in R$, but $x x^2 = x^3 \notin J$.)  
\end{remark}

Despite this fact, there is a special case when the implication does hold, and this will be of use to us.

\begin{lemma} \label{local-units-imply-trans}
Let $R$ be a ring and let $I$ be an ideal of $R$ with the property that $I$ has a set of local units.  If $J$ is an ideal of $I$, then $J$ is an ideal of $R$.
\end{lemma}

\begin{proof}
Let $r \in R$ and $x \in J$.  Since $I$ has a set of local units, there exists $t \in I$ with $tx = x$.  Because $I$ is an ideal, we have that $rt \in I$.  Hence $rx = r(tx) = (rt) x \in J$.  A similar argument shows that $xr \in I$.
\end{proof}

\section{The Graded Uniqueness Theorem} \label{GUT-sec}

In this section we prove a Graded Uniqueness Theorem for Leavitt path algebras.  This result is analogous to the Gauge-Invariant Uniqueness Theorem for graph $C^*$-algebras (\cite[Theorem~2.1]{BPRS} and \cite[Theorem~2.1]{BHRS}).  The main difference is that the grading on the Leavitt path algebra $L_K(E)$ replaces the gauge action of the graph $C^*$-algebra $C^*(E)$.  Consequently, the homogeneous component of degree zero $L_K(E)_0$ replaces the fixed point algebra of the gauge action $C^*(E)^\gamma$, and the hypothesis that our homomorphism is graded replaces the hypothesis that the homomorphism is equivariant with respect to the gauge actions.  

\begin{lemma} \label{ideal-gen-by-zero-part}
Let $I$ be a graded ideal of $L_K(E)$.  Then $I$ is generated as an ideal by the set $I_0 := I \cap L_K(E)_0$. 
\end{lemma}

\begin{proof}
Let $n >0$.  Given $x \in I_n :=I \cap  L_K(E)_n$, we may write $x = \sum_{k=1}^m \alpha_k x_k$ where $x_k \in L_K(E)_0$ for all $k$, $\alpha_k \in E^n$ for all $k$, and $\alpha_i \neq \alpha_j$ for $i \neq j$.  Then for any $1 \leq i \leq m$ we have $$x_i = \alpha_i^* \left(  \sum_{k=1}^m \alpha_k x_k \right) = \alpha_i^* x \in I.$$  Thus $x_i \in I_0$ and $I_n = L_K(E)_n I_0$.  Similarly, $I_{-n} = I_0 L_K(E)_{-n}$.  Since $I$ is a graded ideal, $I = \bigoplus_{n \in \Z} I_n$, and $I$ is generated as an ideal by $I_0$.
\end{proof}

\begin{definition}
For the Leavitt path algebra $L_K(E)$, and for each $n \in \N$ define the following subalgebras of $L_K(E)_0$:
\begin{align*}
\G_n &:= \algspan \{ \alpha \beta^* : \alpha, \beta \in E^n, r(\alpha)=r(\beta) \} \\
\F_n &:=  \algspan \{ \alpha \beta^* : \alpha, \beta \in E^k, r(\alpha)=r(\beta), 0 \leq k \leq n \} 
\end{align*}
\end{definition}

\begin{remark}
It follows from the above definitions that $\F_{n+1} = \G_{n+1} + \F_n$ for each $n \in \N$.  We also see that $\F_n \subseteq \F_{n+1}$ for all $n \in \N$ and $L_K(E)_0 = \bigcup_{n=0}^\infty \F_n$.  Furthermore, $\G_n$ is an ideal of the subalgebra $\F_n$.
\end{remark}

\begin{lemma} \label{G_0-G_1-reg}
For any Leavitt path algebra $L_K(E)$ we have $\G_0 \cap \G_1 = \algspan \{ v : v \in E^0_\reg \}$.  Furthermore, if $\pi : L_K(E) \to A$ is a ring homomorphism with the property that $\pi (v) \neq 0$ for all $v \in E^0$, then $\pi(\G_0) \cap \pi(\G_1) = \pi \left( \algspan \{ v : v \in E^0_\reg \} \right)$.  
\end{lemma}

\begin{proof}
We will show that $\pi(\G_0) \cap \pi(\G_1) = \pi \left( \algspan \{ v : v \in E^0_\reg \} \right)$.  The fact that $\G_0 \cap \G_1 = \algspan \{ v : v \in E^0_\reg \}$ then follows by taking $\pi$ equal to the identity map on $L_K(E)$.

If $v \in E^0_\reg$, then $v \in \G_0$ and $v = \sum_{s(e)=v} ee^*\in \G_1$.  Hence $v \in \G_0 \cap \G_1$ and $\pi \left( \algspan \{ v : v \in E^0_\reg \} \right) \subseteq \pi(\G_0) \cap \pi(\G_1)$.  

Conversely, choose an element $x \in \pi(\G_0) \cap \pi(\G_1)$.  Since $\G^0 = \algspan \{ v : v \in E^0 \}$, we have that $x = \pi(a)$ where $a = \sum_{k=1}^n \lambda_k v_k$ for $\lambda_k \in K$.  Also, since $x \in \pi(\G^1)$, we have $x = \pi \left( \sum_{j=1}^m \mu_j e_j f_j^* \right)$ where $\mu_j \in K$.  Thus $\pi( \sum_{k=1}^n \lambda_k v_k ) = \pi \left( \sum_{j=1}^m \mu_j e_j f_j^* \right)$.  For any $1 \leq k \leq n$, if we multiply each side of this equation on the left by $\pi(v_k)$ we have $\pi(\lambda_k v_k) = \pi \left( \sum_{s(e_j) = v_k} \mu_j e_j f_j^* \right)$.  Because there exist edges with source equal to $v_k$, the vertex $v_k$ is not a sink.  Furthermore, $v_k$ is not an infinite emitter, for if it was, we could find an edge $e \in s^{-1}(v)$ not equal to any of the $e_j$'s, and multiplying each side of the previous equation on the left by $\pi(ee^*)$ would yield $\pi(\lambda_k ee^*) = 0$, which implies that $\pi( r(e) ) = \pi ( \lambda_k^{-1}e^*) \pi (\lambda_k ee^*) \pi(e) = 0$ giving a contradiction.  Therefore $v_k \in E^0_\reg$, and since this is true for all $1 \leq k \leq n$, we have $x = \pi \left( \sum_{k=1}^n \lambda_k v_k \right) \in \pi \left( \algspan \{ v : v \in E^0_\reg \} \right)$.
\end{proof}

\begin{lemma} \label{F_n-replace-G_n-lem}
If $L_K(E)$ is a Leavitt path algebra, then $\F_n \cap \G_{n+1} = \G_n \cap \G_{n+1}$ for all $n \in \N$.  Furthermore, if $\pi : L_K(E) \to A$ is a ring homomorphism, then $\pi( \F_n) \cap \pi(\G_{n+1} )= \pi( \G_n) \cap \pi(\G_{n+1})$ for all $n \in \N$.
\end{lemma}

\begin{proof}
Let $n \in \N$.  We will show that $\pi( \F_n) \cap \pi(\G_{n+1} )= \pi( \G_n) \cap \pi(\G_{n+1})$.  The fact that $\F_n \cap \G_{n+1} = \G_n \cap \G_{n+1}$ then follows by taking $\pi$ equal to the identity map on $L_K(E)$.

Since $\G_n \subseteq \F_n$, we have $\pi(\G_n) \cap \pi(\G_{n+1}) \subseteq \pi(\F_n) \cap \pi(\G_{n+1})$.  Conversely, let $x \in \pi(\F_n) \cap \pi(\G_{n+1})$.  List the elements of $E^n$ as $E^n = \{ \alpha_i : i \in I \}$ where $I = \{1, \ldots, N \}$ if $E^n$ is finite, and $I = \N$ if $E^n$ is infinite.  If we define $t_m := \pi( \sum_{i=1}^m \alpha_i \alpha_i^*)$, then $\{ t_m \}_{m \in I}$ is a set of local units for $\pi(\G_n)$.  Furthermore, since $x \in \pi(\G_{n+1})$ we may write $x = \pi(\sum_{i \in F} \alpha_i x_i)$ for $\alpha_i \in E^n$, $x_i \in L_K(E)$, and some finite set $F \subseteq I$.  Thus there exists $m \in I$ such that $t_m x = x$.  Because $t_m \in \pi(\G_n)$ and $\pi(\G_n)$ is an ideal in $\pi(\F_n)$, it follows that $x = t_m x \in \pi(\G_n)$.  Thus $x \in \pi(\G_n) \cap \pi(\G_{n+1})$.
\end{proof}

\begin{lemma} \label{pi-dist-over-cap-lem}
For any Leavitt path algebra $L_K(E)$, if $\pi : L_K(E) \to A$ is a ring homomorphism with the property that $\pi (v) \neq 0$ for all $v \in E^0$, then $\pi (\F_n \cap \G_{n+1}) = \pi (\F_n) \cap \pi(\G_{n+1})$.
\end{lemma}

\begin{proof}
By Lemma~\ref{F_n-replace-G_n-lem} it suffices to show that $\pi (\G_n \cap \G_{n+1}) = \pi (\G_n) \cap \pi(\G_{n+1})$.  We trivially have that $\pi (\G_n \cap \G_{n+1}) \subseteq \pi (\G_n) \cap \pi(\G_{n+1})$.  

For the converse, let $x \in \pi (\G_n) \cap \pi(\G_{n+1})$.  Since $x \in \pi(\G_n)$ we may write $x = \pi(a)$ where $a = \sum_{k=1}^m \lambda_k \alpha_k \beta_k^*$ for $\lambda_k \in K$ and $\alpha_k, \beta_k \in E^n$.  For any $1 \leq k \leq n$, if we multiply the equation $\pi(a) = \pi \left( \sum_{k=1}^m \lambda_k \alpha_k \beta_k^* \right)$ on the left by $\pi(\alpha_k^*)$ and on the right by $\pi(\beta_k)$ we obtain $\pi( \alpha_k^* a \beta_k) = \pi( \lambda_k r(\alpha_k) ) \in \pi(\G_0)$.  Furthermore, since $x \in \pi(\G_{n+1})$ it follows that $\pi( \alpha_k^* a \beta_k) = \pi(\alpha_k^*)x\pi(\beta_k) \in \pi(\G_1)$.  Therefore, $\pi( \alpha_k^* a \beta_k) \in \pi(\G_0) \cap \pi(\G_1)$.  Because $\pi( \alpha_k^* a \beta_k) = \pi( \lambda_k r(\alpha_k) )$ it follows from Lemma~\ref{G_0-G_1-reg} that $r(\alpha_k) \in E^0_\reg$.  Since this is true for all $k$ we have 
$$a =  \sum_{k=1}^m \lambda_k \alpha_k \beta_k = \sum_{k=1}^m \lambda_k   \alpha_k \left( \sum_{s(e) = r(\alpha_k)} e e^* \right) \beta_k =  \sum_{k=1}^m \sum_{s(e) = r(\alpha_k)} \lambda_k   (\alpha_ke) (\beta_ke)^*$$ so that $a \in \G_{n+1}$.  Hence $a \in \G_n \cap \G_{n+1}$ and $x = \pi(a) \in \pi(\G_n \cap \G_{n+1})$.
\end{proof}

\begin{lemma} \label{F-G-comm-diagram}
Let $n \in \N$ and $\pi : \F_{n+1} \to A$ be a ring homomorphism with the property that $\pi (v) \neq 0$ for all $v \in E^0$.  Also let $\tilde{\pi} : \F_{n+1}/ \G_{n+1} \to \pi(\F_{n+1}) / \pi(\G_{n+1})$ and let $\overline{\pi} : \F_n / (\F_n \cap \G_{n+1}) \to \pi(\F_n) / \pi (\F_n \cap \G_{n+1})$ be the canonical ring homomorphisms induced by $\pi$.  Then there exist isomorphisms $\phi : \F_n / (\F_n \cap \G_{n+1}) \to \F_{n+1} / \G_{n+1}$ and $\phi' :  \pi(\F_n) / \pi(\F_n \cap \G_{n+1}) \to \pi(\F_{n+1}) / \pi(\G_{n+1})$ making the following diagram commute.
\begin{equation*}
\xymatrix{ \F_n / (\F_n \cap \G_{n+1})  \ar[r]^\phi \ar[d]_{\overline{\pi}} & \F_{n+1}/ \G_{n+1} \ar[d]^{\tilde{\pi}} \\
\pi(\F_n) / \pi (\F_n \cap \G_{n+1}) \ar[r]^{\phi'} & \pi(\F_{n+1}) / \pi(\G_{n+1}) \\ }
\end{equation*}
\end{lemma}

\begin{proof}
Define a homomorphism $\phi : \F_n / (\F_n \cap \G_{n+1}) \to \F_{n+1} / \G_{n+1}$ by $\phi(x + (\F_n \cap \G_{n+1})) = x + \G_{n+1}$.  Since $\F_n \cap \G_{n+1} \subseteq \G_{n+1}$, it follows that $\phi$ is well-defined.  In addition, $\phi$ is injective because if $\phi(x + (\F_n \cap \G_{n+1})) = 0 + \G_{n+1}$, then $x \in \G_{n+1}$ and since $x \in \F_n$ also, we have $x \in \F_n \cap \G_{n+1}$ and $x + (\F_n \cap \G_{n+1}) = 0 + (\F_n \cap \G_{n+1})$.  Finally, we see that $\phi$ is surjective because $\F_{n+1} = \G_{n+1} + \F_n$.

We define $\phi' :  \pi(\F_n) / \pi(\F_n \cap \G_{n+1}) \to \pi(\F_{n+1}) / \pi(\G_{n+1})$ by $\phi'(x + \pi(\F_n \cap \G_{n+1})) = x + \pi(\G_{n+1})$.  By Lemma~\ref{pi-dist-over-cap-lem}, we know that $\pi (\F_n \cap \G_{n+1}) = \pi (\F_n) \cap \pi(\G_{n+1})$.  Using this fact, an argument as in the previous paragraph shows that $\phi'$ is an isomorphism.  It is straightforward to check that the diagram in the statement of the lemma commutes.
\end{proof}

\begin{theorem}[Graded Uniqueness Theorem] \label{GUT}
Let $E = (E^0, E^1, r, s)$ be a graph and let $L_K(E)$ be the associated Leavitt path algebra with the usual $\Z$-grading.  If $A$ is a $\Z$-graded ring, and $\pi : L_K(E) \to A$ is a graded ring homomorphism with $\pi(v) \neq 0$ for all $v \in E^0$, then $\pi$ is injective. 
\end{theorem}

\begin{proof}
It follows from Lemma~\ref{ideal-gen-by-zero-part} that the ideal $\ker \pi$ is generated by the set $L_K(E)_0 \cap \ker \pi$.  Thus it suffices to show that the restriction $\pi |_{L_K(E)_0} : L_K(E)_0 \to A$ is injective.  In addition, since $L_K(E)_0 = \bigcup_{n=0}^\infty \F_n$ it suffices to show that the restriction $\pi |_{\F_n} : \F_n \to A$ is injective for all $n \in \N$.  We shall prove this by induction on $n$.

If $n = 0$, then $\F_0 = \algspan \{v : v \in E^0 \}$.  Suppose $\sum_{k=1}^m \lambda_k v_k \in \F_0$ and $\pi(\sum_{k=1}^m \lambda_k v_k) = 0$.  Since the $v_k$'s are mutually orthogonal idempotents, for each $1 \leq j \leq m$ we have $$\pi( \lambda_j v_j) = \sum_{k=1}^m \pi(v_j) \pi( \lambda_k v_k) =  \pi(v_j) \sum_{k=1}^m  \pi(\lambda_k v_k) =  \pi(v_j) \pi \left( \sum_{k=1}^m \lambda_k v_k \right) = 0.$$  Thus $\lambda_j v_j \in I$, and by Lemma~\ref{ring-ideals-alg-ideals-lem} either $\lambda_j = 0$ or $v_j = \lambda_j^{-1} (\lambda_j v_j) \in I$.  Because $\pi(v_j) \neq 0$ by hypothesis, it follows that $v_j \notin I$ and $\lambda_j = 0$.  Because $j$ was arbitrary, we have $\lambda_k = 0$ for all $k$ and $\sum_{k=1}^m \lambda_k v_k = 0$.  Hence $\pi |_{\F_0}$ is injective. 

For the inductive step assume that $\pi |_{\F_n} : \F_n \to A$ is injective.  We then have the following commutative diagram with exact rows
\begin{equation*}
\xymatrix{ 0 \ar[r] & \G_{n+1} \ar[r] \ar[d]^{\pi |_{\G_{n+1}}} & \F_{n+1} \ar[r] \ar[d]^{\pi |_{\F_{n+1}}} & \F_{n+1}/ \G_{n+1} \ar[r] \ar[d]^{\tilde{\pi}} & 0 \\
0 \ar[r] & \pi( \G_{n+1}) \ar[r] & \pi(\F_{n+1}) \ar[r] & \pi(\F_{n+1}) / \pi(\G_{n+1}) \ar[r] & 0 \\ }
\end{equation*}
where $\tilde{\pi} : \F_{n+1} / \G_{n+1} \to \pi(\F_{n+1}) / \pi(\G_{n+1})$ is the canonical homomorphism induced by $\pi$.  We now consider $\pi |_{\G_{n+1}}$.  For each $v \in E^0$ let $$\G_{n+1}(v) := \algspan \{ \alpha \beta^* : \alpha, \beta \in E^{n+1} \text{ and } r(\alpha) = r(\beta) = v \}.$$  Then $\G_{n+1}(v)$ is orthogonal to $\G_{n+1}(w)$ for $v \neq w$, and $\G_{n+1} = \bigoplus_{v \in E^0} \G_{n+1}(v)$ as rings.  Furthermore, for any $\alpha \beta^*, \gamma \delta^* \in  \G_{n+1}(v)$ we have $$\alpha \beta^* \gamma \delta^* = \begin{cases} \alpha \delta^* & \text{ if $\beta = \gamma$} \\ 0 & \text{otherwise.}\end{cases}$$  Hence $\{ \alpha \beta^* : \alpha, \beta \in E^{n+1} \text{ and } r(\alpha) = r(\beta) = v \}$ is a set of matrix units, and $\G_{n+1} \cong M_{m(v)} (K)$ where $m(v)$ is the (possibly infinite) value $m(v) := |\{ \alpha \in E^{n+1} : r(\alpha) = v \} |$.  It follows that $\G_{n+1}$ is simple.  If we let $I = \ker \pi |_{\G_{n+1}}$, then $I = \bigoplus_{v \in E^0} I \cap \G_{n+1}(v)$ (since $\G_{n+1} = \bigoplus_{v \in E^0} \G_{n+1}(v)$ as rings), and by the simplicity of $\G_{n+1}(v)$ the ideal $I \cap \G_{n+1}(v)$ is either $\{ 0 \}$ or all of $\G_{n+1}(v)$.  Furthermore, for each $v \in E^0$, we see that if $\alpha \in E^*$ with $r(\alpha) = v$, then $\pi(\alpha^*\alpha) = \pi(v) \neq 0$ implies $\pi(\alpha \alpha^*) \neq 0$.  Thus $\alpha \alpha^* \notin I \cap \G_{n+1}(v)$, and $I \cap \G_{n+1}(v) = \{ 0 \}$ for all $v \in E^0$.  Hence $I = \{0 \}$ and $\pi |_{\G_{n+1}}$ is injective.

In addition, let $\overline{\pi} : \F_n / (\F_n \cap \G_{n+1}) \to \pi(\F_n) / \pi(\F_n \cap \G_{n+1})$ be the canonical homomorphism induced by $\pi$.  We showed in the previous paragraph that $\pi |_{\G_{n+1}}$ is injective.  Hence $\pi |_{\F_n \cap \G_{n+1}} : \F_n \cap \G_{n+1} \to \pi (\F_n \cap \G_{n+1})$ is an isomorphism.  Since $\pi |_{\F_n}$ is injective by hypothesis, it follows that $\overline{\pi} : \F_n / (\F_n \cap \G_{n+1}) \to \pi(\F_n) / \pi(\F_n \cap \G_{n+1})$ is injective. Therefore Lemma~\ref{F-G-comm-diagram} implies that $\tilde{\pi}$ is injective.

Since $\pi |_{\G_{n+1}}$ and $\tilde{\pi}$ are injective, the commutative diagram above together an application of the Five Lemma shows that $\pi |_{\F_{n+1}} : \F_{n+1} \to \pi(\F_{n+1})$ is injective.  Hence by the Principle of Mathematical Induction, $\pi |_{\F_n} : \F_n \to A$ is injective for all $n \in \N$, and $\pi |_{L_K(E)_0} : L(K)_0 \to A$ is injective.
\end{proof}

\begin{remark}
Note that in Theorem~\ref{GUT} we assumed $\pi$ was a ring homomorphism and not an algebra homomorphism.
\end{remark}

\section{Graded Ideals of $L_K(E)$} \label{graded-ideals-sec}

We shall use the Graded Uniqueness Theorem to characterize the graded ideals of $L_K(E)$.  A characterization of the graded ideals of $L_K(E)$ when $E$ is a row-finite graph was obtained in \cite[Theorem~5.3]{AMP2}, \cite[Lemma~2.3]{PinParMol}, and \cite[Lemma~2.4]{PinParMol}.  As we shall see, in analogy with graph $C^*$-algebras, the description of the graded ideals for Leavitt path algebras of non-row-finite graphs will be more complicated than in the row-finite case --- we will need to use not only saturated hereditary subsets $H$ of vertices, but admissible pairs $(H,S)$ of vertices.

Furthermore, we mention that our method of proof involves a straightforward application of the Graded Uniqueness Theorem.  This is unlike the proof of \cite[Theorem~5.3]{AMP2}, which uses $K$-theory and the order ideals of a monoid $M_E \cong V(L_K(E))$.  Thus in the special case of row-finite graphs, our techniques yield a simpler method for characterizing the graded ideals of $L_K(E)$ than that found in the proof of \cite[Theorem~5.3]{AMP2}.

\begin{definition}
A subset $H \subseteq E^0$ is said to be \emph{hereditary}
if for any $e \in E^1$ we have that $s(e) \in H$ implies
$r(e) \in H$.  A hereditary subset $H \subseteq E^0$ is
called \emph{saturated} if whenever $ 0 < | s^{-1}(v) | <
\infty$, then $\{ r(e) \in H : e \in E^1 \text{ and } s(e)=v
\} \subseteq H$ implies $v \in H$.
\end{definition}

\begin{definition}
If $H$ is a hereditary set, then we define the \emph{saturation of $H$} to be the smallest saturated hereditary subset $\overline{H}$ that contains $H$.
\end{definition}

\begin{definition}
If $H$ is a hereditary subset we define the \emph{breaking vertices of $H$} to be the set $$B_H := \{ v \in E^0
\backslash H : v \in E^0_\infin \text{ and } 0 < | s^{-1}(e) \cap r^{-1}(E^0 \setminus H)| < \infty 
\}$$ and for any $v \in B_H$ we let $$v^H := v -
\sum_{{s(e)=v} \atop {r(e) \notin H}} ee^*.$$
Note that by the definition of $B_H$ the sum appearing above is finite.  Also note that $v^H \in L_K(E)_0$.
\end{definition}

\begin{definition}
We call $(H,S)$ an \emph{admissible pair} if $H$ is a saturated hereditary subset of $E^0$ and $S \subseteq B_H$.  We let $\La_E$ denote the set of admissible pairs of $E$, and we order these elements by $(H,S) \leq (H',S')$ if and only if $H \subseteq H'$ and $S \subseteq H' \cup S'$.  It turns out the ordered set $\La_E$ is actually a lattice, with upper and lower bounds given by:
$$
(H_1,S_1) \wedge (H_2,S_2) := 
\left((H_1 \cap H_2), \left((S_1 \cup H_1) \cup
(S_2 \cup H_2)\right) \cap
B_{H_1 \cap H_2} \right)
$$
\begin{align*}
(H_1,S_1) \vee (H_2,S_2) := 
\Big(
\left(\overline{H_1 \cup H_2}\right) \cup &\left((S_1 \cup S_2) \cap
B^c_{\overline{H_1 
\cap H_2}}\right), \qquad \\
& \qquad \qquad \qquad 
(S_1 \cup S_2) \cap B_{\overline{H_1 \cup H_2}}
\Big).
\end{align*}

\end{definition}

\begin{definition}
If $H$ is a saturated hereditary subset of $E^0$ and $S
\subseteq B_H$, let $I_{(H,S)}$ denote the ideal in $L_K(E)$
generated by $\{ v : v \in H \} \cup \{ v^H : v \in
S\}$.  
\end{definition}

\begin{lemma} \label{lem-ideal-span}
If $H$ is a saturated hereditary subset of $E^0$ and $S
\subseteq B_H$, then $$I_{(H,S)} = \algspan
\big( \{ \alpha \beta^* : r(\alpha)=r(\beta) \in H \}
\cup \{ \alpha v^H \beta^* : r(\alpha) = r(\beta) = v
\in S \} \big)$$ and $I_{(H,S)}$ is a graded ideal of $L_K(E)$ that is self-adjoint.  Moreover, the ideal $I_{(H,S)}$ is an idempotent ring.
\end{lemma}

\begin{proof}
Let $J$ denote the right-hand side of the above equation. 
Since $I_{(H,S)}$ contains $v$ for $v \in H$ and $v^H$
for $v \in S$, we see that $J \subseteq I_{(H,S)}$. 
Conversely, from an examination of the various possibilities
we see that any product of an element in $L_K(E)$ with an element of the form
$\alpha \beta^*$ or $\gamma v^H \delta^*$ will
again be of one of these forms.  Therefore, $J$ is an ideal
which contains the generators of $I_{(H,S)}$ and we deduce
that $J = I_{(H,S)}$.

To see that $I_{(H,S)}$ is graded it suffices to
notice that $\alpha \beta^*$ and
$\alpha v^H \beta^*$ are homogeneous of degree $|\alpha|-|\beta|$.  In addition, we see that $I_{(H,S)}$ is self-adjoint because $\overline{\alpha \beta^*} = \beta \alpha^*$ and  $\overline{\alpha v^H \beta^*} = \beta v^H \alpha^*$.

Finally, to see that $I_{(H,S)}$ is an idempotent ring, simply note that if $\alpha \beta^* \in I_{(H,S)}$ with $r(\alpha) = r(\beta) = v \in H$, then $\alpha \beta^* = (\alpha v) (v \beta^*)$ and we have $\alpha v \in I_{(H,S)}$ and $ v \beta^* \in I_{(H,S)}$.  Likewise, if $\alpha v^H \beta^* \in I_{(H,S)}$ with $r(\alpha) = r(\beta) \in S$, then $\alpha v^H \beta^* = (\alpha v^H) (v^H \beta^*)$ and we have $\alpha v^H \in I_{(H,S)}$ and $v^H \beta^* \in I_{(H,S)}$.  Consequently every element $x \in I_{(H,S)}$ can be written as $x = a_1 b_1 + \ldots + a_n b_n$ for $a_1, \ldots, a_n, b_1, \ldots, b_n \in I_{(H,S)}$, and $I_{(H,S)}$ is an idempotent ring.
\end{proof}

\begin{theorem} \label{ideal-structure-theorem}
Let $E = (E^0,E^1,r,s)$ be a directed graph, and let $L_K(E)$ be the associated Leavitt path algebra.  Let $\La_E$ be the lattice of admissible pairs of $E$, and for $(H,S) \in \La_E$ let $I_{(H,S)}$ denote the ideal generated by $\{ v : v \in H \} \cup \{ v^H : v \in S\}$.  Then
\begin{enumerate}
\item The map $(H,S) \mapsto I_{(H,S)}$ is an isomorphism
from the lattice $\La_E$ onto the lattice of graded ideals of $L_K(E)$. \label{one}
\item For any admissible pair $(H,S) \in \La_E$ we have that $L_K(E) / I_{(H,S)}$ is canonically isomorphic to $L_K(E \setminus (H,S))$, where $E \setminus (H,S)$ is the graph defined by 
\begin{align*}
\qquad (E \setminus (H,S))^0 & := (E^0 \backslash H) \cup \{ v' : v \in B_H \setminus S \} \\ 
\qquad (E \setminus (H,S))^1 & := \{e \in E^1 : r(e) \notin H \} \cup \{e' :
e \in E^1, r(e) \in B_H \backslash S \}
\end{align*}
and $r$ and $s$ are extended to $(E \setminus (H,S))^0$ by setting $s(e') = s(e)$ and $r(e') =
r(e)'$. \label{two}
\item  For any admissible pair $(H,S) \in \La_E$ the ideal $I_{(H,S)}$ and the Leavitt path algebra $L_K (E_{(H,S)})$ are Morita equivalent as rings, where $E_{(H,S)}$ is the graph defined by 
\begin{align*}
\qquad E_{(H,S)}^0 & :=  H \cup S \\ 
\qquad E_{(H,S)}^1 & := \{e \in E^1 : s(e) \in H \} \cup \{e \in E^1 : s(e) \in S \text{ and } r(e) \in H \}
\end{align*}
and we restrict $r$ and $s$ to $E_{(H,S)}^1$. \label{three}
\item  For any hereditary subset $X \subseteq E^0$, if we let $I_X$ denote the ideal in $L_K(E)$ generated by $\{ v : v \in X \}$, then $I_X = I_{(\overline{X},\emptyset)}$.  Moreover, $I_X$ and $L_K (E_X)$ are Morita equivalent as rings, where $E_X := E_{(X, \emptyset)}$ is the graph defined in (\ref{three}). \label{four}
\end{enumerate}
\end{theorem}

\begin{proof}
We shall begin by showing that the
set $\{ v \in E^0 : v \in I_{(H,S)} \}$ is precisely $H$
and the set $\{ v \in B_H : v^H \in I_{(H,S)} \}$ is
precisely $S$.  To begin, we trivially have that $H
\subseteq \{ v \in E^0 : v \in I_{(H,S)} \}$ and $S
\subseteq \{ v \in B_H : v^H \in I_{(H,S)} \}$.  For the
reverse inclusion, let $E \setminus (H,S)$ be the graph of
(\ref{two}) and let $\{a_v : v 
\in (E \setminus (H,S))^0 \} \cup \{b_e, b_{e^*} : e \in (E \setminus (H,S))^1\}$ be the elements of 
$L_K(E \setminus (H,S))$ defined by $$a_v :=
\begin{cases} v & \text{if $v \in (E^0 \backslash H)
\backslash (B_H \backslash S)$} \\
v +v' & \text{if $v \in B_H \backslash S$} \\
0 & \text{if $v \in H$}, \end{cases}$$ $$b_e :=
\begin{cases} e & \text{if $r(e) \in (E^0 \backslash H)
\backslash (B_H \backslash S)$} \\ 
e + e' & \text{if $r(e) \in B_H \backslash S$} \\ 0 &
\text{if $r(e) \in H$} \end{cases}$$ and $$b_{e^*} :=
\begin{cases} e^* & \text{if $r(e) \in (E^0 \backslash H)
\backslash (B_H \backslash S)$} \\ 
e^* + (e')^* & \text{if $r(e) \in B_H \backslash S$} \\ 0 &
\text{if $r(e) \in H$}. \end{cases}$$Then the elements $\{a_v, b_e, b_{e^*} \}$ satisfy the Leavitt path algebra relations for $E$: to see this, we need to use the fact that $H$ is
hereditary to get the Cuntz-Krieger relations at vertices in
$H$, and that $H$ is saturated to see that there are no vertices in
$E^0 \backslash H$ at which a new Cuntz-Krieger relation is
being imposed (in other words, that any singular vertex in
$E \setminus (H,S)^0$ corresponds to a singular vertex in $E^0$).  The
universal property of $L_K(E)$ then gives a homomorphism
$\pi : L_K(E) \rightarrow L_K(E \setminus (H,S))$ with $\pi(v) = a_v$, $\pi(e) = b_e$, and $\pi(e^*) = b_{e^*}$.  Then $\pi$ vanishes on $I_{(H,S)}$ because it kills all the generators $\{ v : v
\in H \} \cup \{ v^H : v \in S\}$.  But $\pi(v) = a_v
\neq 0$ for $v \notin H$, so $v \notin H$ implies $v
\notin I_{(H,S)}$.  Thus $\{ v \in E^0 : v \in I_{(H,S)}
\} \subseteq H$.  Likewise, if $v \in B_H \backslash S$,
then $\pi (v^H) = v' \neq 0$ and $\{ v \in B_H : v^H
\in I_{(H,S)} \} \subseteq S$.

\smallskip

\noindent \underline{\textsc{Proof of (\ref{two}):}}  We shall show that $L_K(E) / I_{(H,S)} \cong
L_K(E \setminus (H,S))$.  Let $\{v \in E^0 \} \cup \{e \in E^1\}$ be the generators for $L_K(E)$.
For each $v \in (E \setminus (H,S))^0$ and $e \in (E \setminus (H,S))^1$ define $$A_v :=
\begin{cases} \qquad \qquad v & \text{if $v \in (E^0 \backslash H) \backslash (B_H \backslash S) $} \\
\sum\limits_{ \{e \in (E \setminus (H,S))^1 : s(e) = v \} } e e^* &
\text{if $v \in B_H
\backslash S$} \\ 
\qquad \qquad v^H & \text{if $v=v'$.}
\end{cases}$$ and
$$B_e := \begin{cases} e & \text{if $r(e) \notin
B_H \backslash S$} \\ 
eA_{r(e)} & \text{if $r(e) \in B_H \backslash S$} \\
eA_{r(e)'} & \text{if $e=e'$.} 
\end{cases}$$ and 
$$B_{e^*} := \begin{cases} e^* & \text{if $r(e) \notin
B_H \backslash S$} \\ 
A_{r(e)}e^* & \text{if $r(e) \in B_H \backslash S$} \\
A_{r(e)'}e^* & \text{if $e=e'$.} 
\end{cases}$$
One can verify that $\{ A_v  + I_{(H,S)} , B_e + I_{(H,S)}, A_v  + I_{(H,S)}, B_e + I_{(H,S)}
\}$ is a set of elements in $L_K(E) / I_{(H,S)}$ satisfying the Leavitt path algebra relations for $E \setminus (H,S)$.  Thus there exists a homomorphism $\phi : L_K(E \setminus (H,S)) \to L_K(E) / I_{(H,S)}$ taking
the generators of $L_K(E \setminus (H,S))$ to the corresponding
elements of $\{ A_v  + I_{(H,S)}, B_e + I_{(H,S)}, B_{e^*} + I_{(H,S)}  \}$.  It
follows from the first paragraph of this proof that $q_v \notin I_{(H,S)}$ for all $v \in (E \setminus (H,S))^0$ and thus the elements $\{ q_v  + I_{(H,S)} \}$ are all
nonzero in $L_K(E) / I_{(H,S)}$.  Furthermore, since
$I_{(H,S)}$ is a graded ideal, it follows that the quotient $L_K(E) / I_{(H,S)}$ is graded.  Since
$\phi$ respects this grading, $\phi$ is a graded homomorphism.   It follows from
Theorem~\ref{GUT} that $\phi$ is injective.  Finally, we
observe that $L_K(E) / I_{(H,S)}$ is generated by $\{ e +
I_{(H,S)} : r(e) \notin H \} \cup \{ v + I_{(H,S)} : v
\notin H \}$.  But for these elements $$v = \begin{cases}
A_v & \text{if $v \notin B_H \backslash S $} \\ A_v + A_{v'}
& \text{if $v \in B_H \backslash S$} \end{cases} \qquad e = \begin{cases} B_e & \text{if $r(e) \notin B_H
\backslash S $} \\ B_e + B_{e'} & \text{if $r(e) \in B_H
\backslash S$} \end{cases}$$ and
$$ e^* = \begin{cases} B_{e^*} & \text{if $r(e) \notin B_H
\backslash S $} \\ B_{e^*} + B_{(e')^*} & \text{if $r(e) \in B_H
\backslash S$} \end{cases}$$
and thus $\phi$ is surjective. 
Hence $L_K(E) / I_{(H,S)} \cong L_K(E \setminus (H,S))$.

\smallskip

\noindent \underline{\textsc{Proof of (\ref{one}):}}  We shall show that $(H,S) \mapsto I_{(H,S)}$ is a lattice
isomorphism.  To see that it is surjective let $I$ be a
graded ideal in $L_K(E)$, and set $H := \{v \in E^0
: v \in I\}$ and $S:= \{v \in B_H : v^H \in I\}$.  Since
$I_{(H,S)} \subseteq I$, we see that $I_{(H,S)}$ and $I$
contain the same $v$'s and $v^H$'s.  Therefore, just as
in the proof of Part (\ref{two}), we see that $L_K(E) /
I_{(H,S)}$ and $L_K(E) / I$ are generated by nonzero elements satisfying the Cuntz-Krieger relations for 
$E \setminus (H,S)$.  Since both $I_{(H,S)}$ and $I$ are graded, both
quotients are graded.  Thus Theorem~\ref{GUT}
implies that the quotient map $\pi : L_K(E \setminus (H,S)) \cong
L_K(E) / I_H \rightarrow L_K(E) / I$ is an isomorphism. 
Hence $I = I_{(H,S)}$.

The fact that $(H,S) \mapsto I_{(H,S)}$ is injective
follows immediately from the fact we deduced in the first paragraph of this proof: the
set $\{ v \in E^0 : p_v \in I_{(H,S)} \}$ is precisely $H$
and the set $\{ v \in B_H : p_v^H \in I_{(H,S)} \}$ is
precisely $S$.  Thus the correspondence $(H,S) \mapsto I_{(H,S)}$ is bijective, and
since it also preserves containment it is a lattice
isomorphism.

\smallskip 

\noindent \underline{\textsc{Proof of (\ref{three}):}}  List the elements of $H = \{ v_1, v_2, \ldots \}$, and list the elements of $S = \{ w_1, w_2, \ldots \}$.  (Each of these sets may be finite or infinite.)  For $i \in \N$ define $$t_i := \begin{cases} v_i & \text{ if $i \leq |H|$} \\ 0 & \text{ if $i > |H|$} \end{cases} \qquad \text{ and } \qquad u_i := \begin{cases} w_i^H & \text{ if $i \leq |S|$} \\ 0 & \text{ if $i > |S|$.} \end{cases} $$  
Notice that the $t_i$'s (respectively, the $u_i$'s) will contain zero terms if and only if $H$ (respectively, $S$) is finite.  Consider the ascending family of idempotents $\{ e_n \}_{n=1}^\infty$, where $e_n := \sum_{i=1}^n t_i + \sum_{i=1}^n u_i$.

If we consider the elements $\{ v : v \in H \} \cup \{v^H : v \in S\}$ and $\{ e, e^* : e \in E^1 \text{ and } s(e) \in H\} \cup \{e, e^* : e \in E^1, s(e) \in S \text{ and } r(e) \in H \}$ in $L_K(E)$, we see that they satisfy the Leavitt path algebra relations for $E_{(H,S)}$ and thus there exists a homomorphism $\pi : L_K (E_{(H,S)}) \to L_K(E)$ taking the generators of $L_K(E_{(H,S)})$ to these elements.  Since this homomorphism is graded, Theorem~\ref{GUT} shows that $\pi$ is injective.  Hence we may identify $L_K(E_{(H,S)})$ with the subalgebra $$\algspan \left( \{ \alpha \beta^* : \alpha, \beta \in E_{(H,S)}^* \text{ and } r(\alpha)=r(\beta) \in H  \} \cup \{ v^H : v \in S \} \right) $$ of $L_K(E)$.

With this identification, we see that $L_K( E_{(H,S)}) = \sum_{n=1}^\infty e_n L_K(E) e_n$.  Moreover, Lemma~\ref{lem-ideal-span} shows that $I_{(H,S)} = \sum_{n=1}^\infty L_K(E) e_n L_K(E)$.  Consider $$\left(\sum_{n=1}^\infty e_n L_K(E) e_n, \sum_{n=1}^\infty L_K(E) e_n L_K(E), \sum_{n=1}^\infty L_K(E)e_n, \sum_{n=1}^\infty e_n L_K(E) , \psi, \phi \right)$$ where $\psi ( m \otimes n) = mn$ and $\phi (n \otimes m) = nm$.  It is straightforward to show that this is a (surjective) Morita context for the idempotent rings $L_K(E_{(H,S)})$ and $I_{(H,S)}$.  It then follows from \cite[Proposition~2.5]{GS} and \cite[Proposition~2.7]{GS} that $L_K(E_{(H,S)})$ and $I_{(H,S)}$ are Morita equivalent.

\smallskip

\noindent \underline{\textsc{Proof of (\ref{four}):}}
We first show that $I_X = I_{(\overline{X}, \emptyset)}$.  As in the proof of \cite[Lemma~4.2]{BPRS}, we note that $\{ v : v \in I_X \}$ is a saturated hereditary subset containing $X$.  Therefore $\overline{X} \subseteq \{ v : v \in I_X \}$, and $I_{(\overline{X},\emptyset)} = \algspan \{ \alpha \beta^* : r(\alpha)=r(\beta) \in \overline{X} \} \subseteq I_X$ by Lemma~\ref{lem-ideal-span}.  But since $I_{(\overline{X},\emptyset)}$ is an ideal containing $\{ v : v \in X \}$, we have that $I_X = I_{(\overline{X},\emptyset)}$.

To see that $I_X$ is Morita equivalent to $L_K(E_X)$, note that as in the proof of (\ref{three}), we may identify $L_K(E_X)$ with the subalgebra of $L_K(E)$ generated by $\{ v : v \in X \} \cup \{ e \in E^1 : s(e) \in X\}$.  If we list the elements of $X = \{ v_ 1, v_2, \ldots \}$, let $$\Lambda := \begin{cases} \{1, 2, \ldots, |X| \} & \text{ if $X$ is finite} \\ \{ 1, 2, \ldots \} & \text{ if $X$ is infinite} \end{cases}$$
and let $e_n := \sum_{i=1}^n v_i$, then  we see that $L_K(E_X) = \sum_{n \in \Lambda} e_n L_K(E) e_n$ and $I_X = \sum_{n \in \Lambda} L_K(E) e_n L_K(E)$.  In addition,
$$\left(\sum_{n \in \Lambda} e_n L_K(E) e_n, \sum_{n \in \Lambda} L_K(E) e_n L_K(E), \sum_{n \in \Lambda} L_K(E)e_n, \sum_{n \in \Lambda} e_n L_K(E), \psi, \phi \right)$$ with $\psi(m \otimes n) = mn$ and $\phi(n \otimes m) =nm$ is a (surjective) Morita context for the idempotent rings $L_K(E_X)$ and $I_X$.  It then follows from \cite[Proposition~2.5]{GS} and \cite[Proposition~2.7]{GS} that $L_K(E_X)$ and $I_X$ are Morita equivalent.

\end{proof}

\begin{example}
Let $E$ be the graph

$$
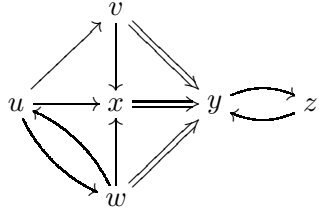
\xymatrix{
 & v \ar[d] \ar@{=>}[dr] & & \\
u \ar[ru] \ar[r] \ar@/_/[rd] & x \ar@{=>}[r]  & y \ar@/^/[r] & z \ar@/^/[l] \\
 & w \ar[u] \ar@/_/[lu] \ar@{=>}[ur] & & \\
}
$$
where the double arrows $\Longrightarrow$ indicate that there are a countably infinite number of edges from $v$ to $y$, from $x$ to $y$, and from $w$ to $y$.  Then $H := \{ y, z \}$ is a saturated hereditary subset.  We see that $B_H = \{ v, w \}$.  If we let $S: = \{ v \}$, then $(H,S)$ is an admissible pair.  Furthermore, $L_K(E) / I_{(H,S)}$ is isomorphic to the $L_K(E\setminus (H,S))$, where $E\setminus (H,S)$ is the graph
$$
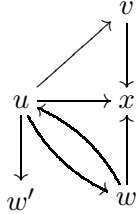
\xymatrix{
 & v \ar[d] \\
u \ar[ru] \ar[r] \ar[d] \ar@/_/[rd] & x \\
w' & w \ar[u] \ar@/_/[lu]  \\
}
$$
and the ideal $I_{(H,S)}$ is Morita equivalent to $L_K(E_{(H,S)})$, where $E_{(H,S)}$ is the graph
$$
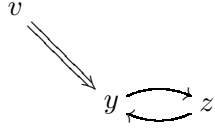
\xymatrix{
v \ar@{=>}[dr] & & \\
 & y \ar@/^/[r] & z \ar@/^/[l] \\
}
$$

\end{example}

\section{The Cuntz-Krieger Uniqueness Theorem}
\label{CK-uniqueness-sec}

In this section we derive another uniqueness theorem for $L_K(E)$, in analogy with the Cuntz-Krieger Uniqueness Theorem for graph $C^*$-algebras.  We show that if the graph $E$ has the property that all closed paths have exits, then we may remove the condition that the homomorphism is graded from the Graded Uniqueness Theorem.  (In other words, when every closed path of $E$ has an exit, a homomorphism on $L_K(E)$ that does not kill any $v$ is injective.)  We then use this to derive a condition on graphs, called Condition~(K), that is equivalent to all ideals in $L_K(E)$ being graded.

\begin{definition}
Let $E$ be a graph.  We sat that a path $\alpha = e_1 \ldots e_n \in E^n$ is a \emph{closed path} if $r(\alpha) = s(\alpha)$, and we say that $\alpha$ is \emph{based at $v$} if $s(\alpha) = v$. We say that an edge $f \in E^1$ is an \emph{exit} for $\alpha$ if $s(f) = s(e_i)$ but $f \neq e_i$ for some $i \in \{1, 2, \ldots, n\}$.
\end{definition}

\begin{definition}
We say that a closed path $\alpha = e_1 \ldots e_n \in E^n$ is \emph{simple} if $s(e_i) \neq s(e_1)$ for $i = 2, 3, \ldots, n$.
\end{definition}

\begin{definition}
A graph $E$ satisfies \emph{Condition~(L)} if every closed path in $E$ has an exit.
\end{definition}

\begin{remark} 
Note that a graph $E$ satisfies Condition~(L) if and only if every closed simple path in $E$ has an exit.
\end{remark}

\begin{lemma} \label{ghost-edges-lem}
Let $E$ be a row-finite graph with no sinks.  If $I$ is a nonzero ideal in $L_K(E)$, then there exists a nonzero element $x \in I$ with the property that $x$ is a polynomial in only ghost edges.
\end{lemma}

\begin{proof}
Choose a nonzero element $y \in I$.  Since the elements of $E^0$ are a set of local units for $L_K(E)$, there exists $v \in E^0$ such that $vy \neq 0$.  We may, by collecting terms, write $$vy = \sum_{\alpha \in F} \alpha z_\alpha$$ for a finite set $F \subseteq E^*$ and nonzero polynomials $z_\alpha$ in only ghost edges.  Let $n = \max \{ | \alpha | : \alpha \in F \}$.

We claim that there exits $\beta \in E^n$ such that $\beta^* vy \neq 0$.  If not, then $\beta^* vy = 0$ for all $\beta \in E^n$ and $\beta \beta^* vy = 0$ for all $v \in E^n$.  Since $E$ is row-finite with no sinks, repeated applications of relation (CK2) shows that $\displaystyle v = \sum_{\{\beta \in E^n: s(\beta)=v\}} \beta\beta^*$.  Thus $$vy = v (vy) = \sum_{\{ \beta \in E^n : s(\beta)=v \} } \beta \beta^* vy = 0,$$ which is a contradiction.  Hence there exists $\beta \in E^n$ such that $\beta^* vy \neq 0$.

Let $x = \beta^* vy$.  Then $x$ is a nonzero element of $I$.  In addition, since $|\beta| \geq |\alpha|$ for all $\alpha \in F$, we see that $\beta^* \alpha$ is a path in ghost edges for all $\alpha \in F$, and $$x = \beta^* vy = \sum_{\alpha \in F} (\beta^* \alpha) z_\alpha$$ is a polynomial in only ghost edges.
\end{proof}

\begin{remark}
If $E$ is a graph, then a desingularization of $E$ is a graph $F$ that is row-finite and has no sinks.  The desingularization was introduced in \cite[\S2]{DT1}, and it was proven in \cite[Theorem~2.11]{DT1} that forming the desingularization preserves the Morita equivalence class of the associated graph $C^*$-algebra.  Abrams and Aranda-Pino have recently shown in \cite[Theorem~5.2]{AbrPino3} that forming the desingularization also preserves the Morita equivalence class of the associated Leavitt path algebra.  We give another proof of this fact, as well as establish an isomorphism between the ideals in the two Leavitt path algebras.
\end{remark}

\begin{lemma} \label{ideal-iso-lem}
Let $E$ be a graph and let $F$ be a desingularization of $E$.  List the vertices of $E$ as $E^0 := \{v_i \}_{i \in \Lambda}$, where $\Lambda := \{1, 2, \ldots, |E^0| \}$ if $E^0$ is finite and $\Lambda := \{1, 2, \dots \}$ if $E^0$ is infinite, and define $e_n := \sum_{i=1}^n v_i$.  Then $L_K(E)$ is isomorphic to the subalgebra $\sum_{n \in \Lambda} e_n L_K(F) e_n$ of $L_K(F)$, and $L_K(E)$ and $L_K(F)$ are Morita equivalent.   Furthermore, the map $I \mapsto \sum_{n \in \Lambda} e_nIe_n$ is a lattice isomorphism from the lattice of ideals of $L_K(F)$ onto the lattice of ideals of $L_K(E)$.
\end{lemma}

\begin{proof}
Let $F$ be a desingularization of $E$, as described in \cite[\S2]{DT1} and \cite[\S5]{AbrPino3}, and assume the reader is familiar with this construction.  As described in the second paragraph of the proof of \cite[Theorem~5.2]{AbrPino3} there is a monomorphism $\phi : L_K(E) \to L_K(F)$ that takes
\begin{align*}
v &\mapsto v \qquad \qquad  \qquad \text{ if $v \in E^0$ }  \\
e &\mapsto e \qquad \qquad  \qquad \text{ if $s(e) \in E^0_\textnormal{reg}$} \\
e_i &\mapsto f_1 \ldots f_{i-1} g_i \qquad \text{ if $e_i \in E^1$ and $s(e_i) \in E^0_\textnormal{sing}$}  \\
e^* &\mapsto e^* \qquad \qquad  \qquad \text{ if $e \in E^1$ and $s(e) \in E^0_\textnormal{reg}$}  \\
e_i &\mapsto g_i^* f_{i-1}^* \ldots f_1^*  \qquad \text{ if $e_i \in E^1$ and $s(e_i) \in E^0_\textnormal{sing}$.} 
\end{align*}
It follows that the image of this monomorphism is equal to $\sum_{n \in \Lambda} e_n L_K(F) e_n$, and thus 
$L_K(E)$ is isomorphic to the subalgebra $\sum_{n \in \Lambda} e_n L_K(F) e_n$ of $L_K(F)$.  We shall identify $L_K(E)$ with this subalgebra.  We then see that $$\left(\sum_{n \in \Lambda} e_n L_K(F) e_n, \sum_{n \in \Lambda} L_K(F) e_n L_K(F), \sum_{n \in \Lambda} L_K(F)e_n, \sum_{n \in \Lambda} e_n L_K(F), \psi, \phi \right),$$ where $\psi (m \otimes n) =mn$ and $\phi(n \otimes m) = nm$, is a (surjective) Morita context for the idempotent rings $L_K(F) = L_K(F) e_n L_K(F)$ and $e_n L_K(F)e_n \cong L_K(E)$.  Thus $L_K(F)$ and $L_K(E)$ are Morita equivalent.  

Furthermore, since $L_K(F)$ and $L_K(E)$ are rings with sets of local units, it follows from Proposition~\ref{ideal-corresp-Mor-eq} and Remark~\ref{ideal-corresp-rings-loc-units} that the map that sends $I$ to $$\phi \left( \sum_{n \in \Lambda} e_n L_K(F) I, \sum_{n \in \Lambda} L_K(F)e_n \right) = \sum_{n \in \Lambda} e_n L_K(F) I \sum_{n \in \Lambda} L_K(F)e_n = \sum_{n \in \Lambda} e_n I e_n$$ is a lattice isomorphism from the lattice of ideals of $L_K(F)$ onto the lattice of ideals of $L_K(E)$.
\end{proof}

\begin{theorem}[Cuntz-Krieger Uniqueness] \label{CK-Uniqueness}
Let $E$ be a graph that satisfies Condition~(L), and let $L_K(E)$ be the associated Leavitt path algebra.  If $\pi : L_K(E) \to A$ is a ring homomorphism with $\pi(v) \neq 0$ for all $v \in E^0$, then $\pi$ is injective.
\end{theorem}

\begin{proof}
First consider the case when $E$ is row-finite with no sinks.  Let $I := \ker \pi$.  If $I$ is nonzero, then by Lemma~\ref{ghost-edges-lem} there is a nonzero element $x \in I$ such that $x$ is a polynomial in only ghost edges.  By \cite[Corollary~3.8]{AbrPino} there exist $v \in E^0$ with $v \in I$.  But this contradicts the assumption that $\pi(v) \neq 0$.  Hence we must have $I = \{ 0 \}$, and $\pi$ is injective.

Now consider the case when $E$ is not necessarily row-finite with no sinks.  Again, let $I := \ker \pi$.  Then $\pi(v) \neq 0$ for all $v \in E^0$ we have that $I \cap E^0 = \emptyset$.  Let $F$ be a desingularization of $E$, list the vertices of $E$ as $E^0 := \{v_i \}_{i \in \Lambda}$, where $$\Lambda := \begin{cases} \{1, 2, \ldots, |E^0| \} & \text{ if $E^0$ is finite} \\ \{1, 2, \ldots \} & \text{ if $E^0$ is infinite} \end{cases}$$ and define $e_n := \sum_{i=1}^n v_i$.  By Lemma~\ref{ideal-iso-lem}, there is an ideal $J$ of $C^*(E)$ such that $\sum_{n \in \Lambda} e_nJe_n = I$.  

We claim that $J \cap F^0 = \emptyset$ by arguing as follows:  Certainly if $v \in E^0$, then $v \in J$ implies that for large enough $n$ we have $v = e_nve_n \in I$, which contradicts the fact that $I \cap E^0 = \emptyset$.  Hence $J \cap E^0 = \emptyset$.  In addition, if $v \in F^0 \setminus E^0$, then $v$ is on a tail added to some singular vertex of $E^0$, and there exists a path $\alpha \in F^*$ with $s(\alpha) = v$ and $r(\alpha) \in E^0$.  But then $$r(\alpha) = \alpha^* \alpha = \alpha^* s(\alpha) \alpha = \alpha^* v \alpha \in J$$ contradicting the fact that $J \cap E^0 = \emptyset$.  Hence we must have that $J \cap F^0 = \emptyset$.

Since $F$ is a desingularization of $E$, the graph $F$ is row-finite with no sinks and it follows from Lemma~2.7(a) that $F$ satisfies Condition~(L).  Therefore, the fact that $J \cap F^0 = \emptyset$, together with Lemma~\ref{ghost-edges-lem} and \cite[Corollary~3.8]{AbrPino}, implies that $J = \{ 0 \}$.  But then $I = \sum_{n \in \Lambda} v_nJv_n = \{0 \}$ and $\pi$ is injective.

\end{proof}

\begin{remark}
Note that in Theorem~\ref{CK-Uniqueness}, as with the Graded Uniqueness Theorem, we assumed $\pi$ was a ring homomorphism and not an algebra homomorphism.
\end{remark}

\begin{corollary}
Let $E$ be a graph that satisfies Condition~(L), and let $L_K(E)$ be the associated Leavitt path algebra.  If $I$ is a nonzero ideal of $L_K(E)$, then $I \cap E^0 \neq \emptyset$.
\end{corollary}

Using the Cuntz-Krieger Uniqueness Theorem we can characterize those graphs whose associated Leavitt path algebras have the property that all ideals are graded.

\begin{definition}
A graph $E$ satisfies \emph{Condition~(K)} if every vertex in $E^0$ is either the base of no closed path or the base of at least two simple closed paths.
\end{definition}

The following proposition is well-known to graph $C^*$-algebraists.  It has been proven in the row-finite case in \cite[Proposition~1.17]{Tom9} and \cite[Theorem~4.5(2),(3)]{PinParMol}, and in the non-row-finite case the proof is nearly identical.

\begin{proposition} \label{K-implies-quotient-L}
If $E$ is a graph, then $E$ satisfies Condition~(K) if and only if for every admissible pair $(H,S)$ the graph $E \setminus (H,S)$ defined Theorem~\ref{ideal-structure-theorem}(\ref{two}) satisfies Condition~(L).
\end{proposition}

\begin{lemma} \label{closed-path-M-Kx}
If $E$ is the graph consisting of a single simple closed path of length $n$; i.e., 
$$
E^0 = \{v_1, \ldots, v_n\} \quad E^1 = \{e_1, \ldots e_n \} $$ 
$$s(e_i) = v_i \quad \text{ for $1 \leq i \leq n$}$$ 
$$r(e_i) = v_{i+1} \quad \text{ for $1 \leq i < n$} \quad \text{ and } \quad r(e_n) = v_1$$
then $L_K(E) \cong M_n(K) \otimes K[x,x^{-1}]$.
\end{lemma}

\begin{proof}
Let $\{ e_{i,j} \}$ denote the matrix units of $M_n(K)$.  Note that the map which sends
\begin{align*}
v_i &\mapsto e_{i,i} \otimes 1 \quad \text{ for $1 \leq i \leq n$} \\
e_i &\mapsto e_{i,i+1} \otimes x \quad \text{ for $1 \leq i < n$} \\
e_n &\mapsto e_{1,n} \otimes x \\
e_i^* &\mapsto e_{i+1,i} \otimes x^{-1} \quad \text{ for $1 \leq i < n$} \\
e_n^* &\mapsto e_{n,1} \otimes x^{-1} 
\end{align*}
is a homomorphism.  (Simply check that the image elements satisfy the defining relations for the generators of $L_K(E)$, and use the universal property.)  Using the natural grading on $M_n(K) \otimes K[x,x^{-1}]$, we see that this homomorphism is graded, and by Theorem~\ref{GUT} the homomorphism is injective.  It is straightforward to show that the homomorphism is surjective, and thus an isomorphism.
\end{proof}

The ideas in the following lemma were suggested to the author by Enrique Pardo.

\begin{lemma} \label{ideals-tran-lem}
Let $E$ be a graph, and let $H$ be a saturated hereditary subset of $E$.  Also let $I_{(H,\emptyset)}$ be the ideal of $L_K(E)$ defined in Theorem~\ref{GUT}.  Then $I_{(H,\emptyset)}$ is a ring with a set of local units.
\end{lemma}

\begin{proof}
Let $$F(X) := \{ \alpha = e_1 \ldots e_{|\alpha|} \in E^* : r(e_{|\alpha|}) \in H \text{ and } r(e_i) \notin H \text{ for } 1 \leq i \leq |\alpha| -1 \}.$$  List the elements of $F(X) = \{ \alpha_1, \alpha_2, \ldots \}$ and list the elements of $H = \{ v_1, v_2, \ldots \}$.  Define
$$e_i := \begin{cases} v_i & \text{ if $i \leq |H|$} \\ 0 & \text{ if $i > |H|$} \end{cases} \qquad f_i := \begin{cases} \alpha_i \alpha_i^* & \text{ if $i \leq |F(X)|$} \\ 0 & \text{ if $i > |F(X)|$}. \end{cases}$$  Note that the $e_i$'s (respectively, the $f_i$'s) will contain zero terms if and only if $|H|$ (respectively, $|F(X)|$) is infinite.  Clearly these elements are idempotents, and one can also see that these idempotents are mutually orthogonal.  (Note that by definition of $F(X)$, no $\alpha_i$ can extend an $\alpha_j$ for $i \neq j$ and thus $\alpha_i \alpha_i^* \alpha_j \alpha_j = 0$.)  Thus if we define $t_n := \sum_{i=1}^n e_i + \sum_{i=1}^n  f_i$, we see that $\{ t_n \}_{n=1}^\infty$ is a set of commuting idempotents.  From Lemma~\ref{lem-ideal-span} we have that $$I_{(H,\emptyset)} = \algspan \{ \alpha \beta^* : r(\alpha) = r(\beta) \in H \}.$$  Since any $\alpha \in E^*$ with $r(\alpha) \in H$ either has $s(\alpha) \in H$ or may be written as $\alpha = \gamma \delta$ for $\gamma \in F(X)$, we see that $\{ t_n \}_{n=1}^\infty$ is a set of local units for $I_{(H,\emptyset)}$.
\end{proof}

\begin{lemma} \label{not-L-non-graded-ideals}
Let $E$ be a graph that contains a closed path with no exit.  Then $L_K(E)$ contains ideals that are not graded.  In fact, the cardinality of the set of ideals in $L_K(E)$ that are not graded will be at least $\max \{ \aleph_0, |K| \}$.  
\end{lemma}

\begin{proof}
Let $\alpha := e_1 \ldots e_n$ be a closed path with no exists in $E$.  If we let $X : = \{s(e_i) \}_{i=1}^n$, then since $\alpha$ has no exits, $X$ is a hereditary subset of $E^0$.  By Theorem~\ref{ideal-structure-theorem}(\ref{four}) $L_K(E_X)$ is Morita equivalent to the ideal $I_X = I_{(\overline{X}, \emptyset)}$ in $L_K(E)$.  However, $E_X$ is the graph which consists of a single closed path, and thus $L_K(E_X) \cong M_n(K) \otimes K[x,x^{-1}]$ by Lemma~\ref{closed-path-M-Kx}.  Theorem~\ref{ideal-structure-theorem}(\ref{one}) implies that $L_K(E)$ has no proper nontrivial graded ideals.  In addition, $K[x, x^{-1}]$ is a Principal Ideal
Domain (in fact, a Euclidean Domain) and a unital ring.  Thus any two-sided ideal of $M_n(K) \otimes K[x, x^{-1}]$ is of the form $M_n(K) \otimes I$ for an ideal $I$ of $K[x, x^{-1}]$, and since any such $I$ is generated by a nonzero polynomial in the variables $x$ and $x^{-1}$ with coefficients from $K$, there are infinitely many ideals in $M_n(K) \otimes K[x, x^{-1}]$ and the cardinality of the set of these ideals corresponds with $\max \{ \aleph_0, |K| \}$.  It follows that $L_K(E_X)$ contains an equal number of ideals that are not graded.  Because the Morita context described in the proof of Theorem~\ref{ideal-structure-theorem}(\ref{four}) gives a lattice isomorphism from ideals of $L_K(E_X)$ to ideals of $I_X$ that preserves the grading, we may conclude that $I_X$ contains at least this many ideals that are not graded.  Since $I_X = I_{(\overline{X}, \emptyset)}$ has a set of local units by Lemma~\ref{ideals-tran-lem}, it follows from Lemma~\ref{local-units-imply-trans} that ideals of $I_X$ are ideals of $L_K(E)$.  Hence $L_K(E)$ contains at least $\max \{ \aleph_0, |K| \}$ ideals that are not graded.  
\end{proof}

These results together with the Cuntz-Krieger Uniqueness Theorem give us the following theorem, which generalizes \cite[Proposition~3.3]{PinParMol}.

\begin{theorem} \label{K-iff-ideals-graded}
If $E$ is a graph, then $E$ satisfies Condition~(K) if and only if every ideal in $L_K(E)$ is graded.
\end{theorem}

\begin{proof}
Suppose that $E$ satisfies Condition~(K).  If $I$ is an ideal of $L_K(E)$, let $H := \{ v : v \in I \}$ and let $S := \{ v : v^H \in I \}$.  Then $I_{(H,S)} \subseteq I$, and we have a canonical surjection $q : L_K(E) / I_{(H,S)} \to L_K(E) / I$.  By Theorem~\ref{ideal-structure-theorem}(\ref{two}) there exists a canonical isomorphism $\phi : L_K(E \setminus (H,S)) \to L_K(E) / I_{(H,S)}$.  Thus the composition $q \circ \phi : L_K(E \setminus (H,S)) \to L_K(E) / I$ has the property that $q \circ \phi (v) \neq 0$ for all $v \in E^0$.  Since $E$ satisfies Condition~(K), it follows from Proposition~\ref{K-implies-quotient-L} that $E \setminus (H,S)$ satisfies Condition~(L).  Hence we may apply Theorem~\ref{CK-Uniqueness} to conclude that $q \circ \phi$ is injective.  Since $\phi$ is an isomorphism, this implies that $q$ is injective and $I = I_{(H,S)}$.  It then follows from Lemma~\ref{lem-ideal-span} that $I$ is graded.

Conversely, suppose that $E$ does not satisfy Condition~(K).  Then there exists an admissible pair $(H,S)$ such that $E \setminus (H,S)$ does not satisfy Condition~(L).  Thus there exists a closed simple path with no exit in $E \setminus (H,S)$, and by Lemma~\ref{not-L-non-graded-ideals} the algebra $L_K(E \setminus (H,S)) \cong L_K(E)/ I_{(H,S)}$ contains an ideal $I$ that is not graded.  If we let $q : L_K(E) \to L_K(E)/ I_{(H,S)}$, then $q$ is graded, and $q^{-1}(I)$ is an ideal of $L_K(E)$ that is not graded.
\end{proof}

\begin{corollary}
If $E$ satisfies Condition~(K), then the map $(H,S) \mapsto I_{(H,S)}$ is a lattice isomorphism from the lattice of admissible pairs of $E$ onto the lattice of ideals of $L_K(E)$.  
\end{corollary}

The following result characterizes simplicity for Leavitt path algebras.  A special case of this result for Leavitt path algebras of row-finite graphs was obtained in \cite[Theorem~3.11]{AbrPino}.  Our method of proof is different, and in the row-finite case gives a shorter proof than the one in \cite{AbrPino}.

\begin{theorem}
Let $E$ be a graph.  The Leavitt path algebra $L_K(E)$ is simple if and only if $E$ satisfies the following two conditions:
\begin{itemize}
\item[(i)] The only saturated hereditary subsets of $E^0$ are $\emptyset$ and $E^0$, and
\item[(ii)] The graph $E$ satisfies Condition~(L).
\end{itemize}
\end{theorem}

\begin{proof}
Suppose that $L_K(E)$ is simple.  Then the only ideals of $L_K(E)$ are $\{0\}$ and $L_K(E)$, both of which are graded.  By Theorem~\ref{K-iff-ideals-graded} we have that $E$ satisfies Condition~(K).  It then follows from Theorem~\ref{ideal-structure-theorem}(\ref{one}) and the simplicity of $L_K(E)$ that the only saturated hereditary subsets of $E^0$ are $\emptyset$ and $E^0$.  Hence (i) holds.  In addition, since Condition~(K) implies Condition~(L) (cf.~Proposition~\ref{K-implies-quotient-L}) we have that (ii) holds.

Conversely, suppose that (i) and (ii) hold.  We shall show that $E$ satisfies Condition~(K).  Let $v$ be a vertex and let $\alpha = e_1 \ldots e_n$ be a closed simple path based at $v$.  By (ii) we know that $\alpha$ has an exit $f$; i.e. there exists $f \in E^1$ with $s(f) = s(e_i)$ and $f \neq e_i$ for some $i$.  If we let $H$ be the set of vertices in $E^0$ such that there is no path from that vertex to $v$, then $H$ is saturated hereditary.  By (i) we must have either $H = \emptyset$ or $H = E^0$.  Since $v \notin H$ it must be the case that $H = \emptyset$.  Hence for every vertex in $E^0$, there is a path from that vertex to $v$.  Choose a path $\beta \in E^*$ from $r(f)$ to $v$ of minimal length.  Then $e_1 \ldots e_{i-1} f \beta$ is a simple closed path based at $v$ that is distinct from $\alpha$.  Hence $E$ satisfies Condition~(K).  It then follows from Theorem~\ref{ideal-structure-theorem}(\ref{one}) and (i) that $L_K(E)$ is simple.
\end{proof}

Condition (i) and (ii) in the above theorem can be reformulated in a number of equivalent ways.  The equivalence of the statements (2)--(5) in Proposition~\ref{simple-equiv-prop} are elementary facts about directed graphs (cf.~\cite[Theorem~1.23]{Tom9} and \cite[Proposition~3.2]{AbrPino3}).

\begin{definition}
A graph $E$ is \emph{cofinal} if whenever $e_1 e_2 e_3 \ldots$ is an infinite path in $E$ and $v \in E^0$, then there exists a finite path from $v$ to $s(e_i)$ for some $i \in \N$.
\end{definition}

\begin{proposition} \label{simple-equiv-prop}
Let $E$ be a graph, and let $L_K(E)$ be the associated Leavitt path algebra.  Then the following are equivalent.
\begin{enumerate}
\item $L_K(E)$ is simple.
\item $E$ satisfies Condition~(L), and the only saturated hereditary subsets of $E^0$ are $\emptyset$ and $E^0$.
\item $E$ satisfies Condition~(K), and the only saturated hereditary subsets of $E^0$ are $\emptyset$ and $E^0$.
\item $E$ satisfies Condition~(L), $E$ is cofinal, and whenever $v \in E^0_\textnormal{sing}$ and $w \in E^0$ there is a path from $v$ to $w$.
\item $E$ satisfies Condition~(K), $E$ is cofinal, and whenever $v \in E^0_\textnormal{sing}$ and $w \in E^0$ there is a path from $v$ to $w$.
\end{enumerate}
\end{proposition}

\section{The Leavitt path algebra $L_\C(E)$ and graph $C^*$-algebra $C^*(E)$}
\label{Leavitt-and-graph-algebras-sec}

Much of the progress made in the study of Leavitt path algebras has been motivated by the theory developed for graph $C^*$-algebras in the past decade.  Although Leavitt path algebras and graph $C^*$-algebras have much in common (e.g., they each have generators satisfying the same set of relations, and many similar results hold for the objects of each class), there are important differences: 

\begin{enumerate}
\item[(1)] Graph $C^*$-algebras are algebras over $\C$ with analytic structure, while Leavitt path algebras are algebras (without any analytic structure) over an arbitrary field $K$.
\item[(2)] Recently it has been shown that certain graph $C^*$-algebra results do not hold for Leavitt path algebras (e.g., it is shown in \cite{AraPar} that there is a graph $E$ whose associated Leavitt path algebra $L_K(E)$ has stable rank 2 but whose associated graph $C^*$-algebra $C^*(E)$ has stable rank 1).
\item[(3)] When similar results do hold for graph $C^*$-algebras and Leavitt path algebras, there is no obvious way to deduce the results for one class from the results of the other class.  Also, in the current literature it is not uncommon for the proof of a result for one class to be substantially different from the proof of the corresponding result for the other class.
\end{enumerate}

\noindent In this section we examine the relationship between the Leavitt path algebra $L_\C(E)$ and the graph $C^*$-algebra $C^*(E)$ for a fixed graph $E$.

\begin{definition} \label{CK-relations}
Let $E$ be a graph.  A Cuntz-Krieger $E$-family is a collection of mutually orthogonal projections $\{ p_v : v \in E^0 \}$ and a collection of partial isometries with mutually orthogonal ranges $\{ s_e : e \in E^1 \}$ satisfying the following three relations:
\begin{itemize}
\item[(CK1)] $s_e^*s_e = p_{r(e)}$ \text{ for all $e \in E^1$}
\item[(CK2)] $p_v = \sum_{s(e)=v} s_es_e^*$ \text{ for all $v \in E^0_\textnormal{reg}$}
\item[(CK3)] $s_es_e^* \leq p_{s(e)}$ \text{ for all $e \in E^1$}.
\end{itemize}
\end{definition}

\begin{definition}
If $E$ is a graph, the \emph{graph $C^*$-algebra} $C^*(E)$ is the $C^*$-algebra generated by a universal Cuntz-Krieger $E$-family; that is, $C^*(E)$ is generated by a Cuntz-Krieger $E$-family $\{s_e, p_v \}$, and whenever $\{t_e, q_v \}$ is a Cuntz-Krieger $E$-family sitting inside a $C^*$-algebra $A$, then there exists a $*$-homomorphism $\phi : C^*(E) \to A$ with $\phi(s_e) = t_e$ for all $e \in E^1$ and $\phi (p_v) = q_v$  for all $v \in E^0$.
\end{definition}

The existence and uniqueness (up to isomorphism) of the graph $C^*$-algebra is proven in \cite[Proposition~1.21]{Rae}.  If $\alpha = e_1 \ldots e_n \in E^*$ we write $s_\alpha$ for the product $s_{e_1} \ldots s_{e_n}$.

It has frequently been stated (without proof) that if $E$ is a graph, then the Leavitt path algebra $L_\C(E)$ is isomorphic to the dense $C^*$-subalgebra 
\begin{equation} \label{dense-subalg}
\mathcal{A} := \algspan \{ s_\alpha s_\beta^* : \alpha, \beta \in E^* \text{ and } r(\alpha) = r(\beta) \}
\end{equation}
 of $C^*(E)$.  (This is asserted in the second paragraph of the introduction to \cite{AbrPino}, and the second paragraph of the introduction to \cite{AMP2}, among other places.  It is also used implicitly throughout much of the work in \cite{AMP2}.)  In personal communication with the author, Ara and Pardo explained that when they used this fact in \cite{AMP2}, they knew it was true because it is a consequence of their description of graded ideals in Leavitt path algebras of row-finite graphs \cite[Theorem~5.3]{AMP2}.  Thus they deduced the row-finite case of the result as a consequence of the Graded Uniqueness Theorem for row-finite Leavitt path algebras.  We shall use our Graded Uniqueness Theorem, stated in Theorem~\ref{GUT}, to give a proof of the result in the general case.  As far as we know, this is the only place (even for the row-finite case) that a proof has been written down.  Moreover, as we shall see, this has consequences for the construction of the graph $C^*$-algebra.

Let us examine the statement ``the Leavitt path algebra $L_\C(E)$ is naturally isomorphic to the dense $*$-subalgebra $\mathcal{A}$ of $C^*(E)$ described in (\ref{dense-subalg})" and consider the subtleties that make it non-obvious.  To begin, let us carefully state the universal properties of $L_K(E)$ and $C^*(E)$.

\smallskip

\noindent \textbf{Universal property of $\boldsymbol{L_K(E)}$:}  If $A$ is a $K$-algebra containing a collection of elements $\{a_v : v \in E^0 \} \cup \{b_e, b_{e^*} : e \in E^1\}$ satisfying the relations in Definition~\ref{Leavitt-def}, then there is a homomorphism $\phi : L_K(E) \to A$ with $\phi(v) = b_v$ for all $v \in E^0$ and $\phi(e) = b_e$ and $\phi(e^*) = b_{e^*}$ for all $e \in E^1$.

\smallskip

\noindent \textbf{Universal property of $\boldsymbol{C^*(E)}$:}  If $A$ is a $C^*$-algebra containing a collection of projections $\{q_v : v \in E^0 \}$ and partial isometries with mutually orthogonal ranges $\{t_e : e \in E^1\}$ satisfying the relations in Definition~\ref{CK-relations}, then there is a $*$-homomorphism $\phi : C^*(E) \to A$ with $\phi(p_v) = q_v$ for all $v \in E^0$ and $\phi(s_e) = t_e$ for all $e \in E^1$.

\smallskip

If we look at $\mathcal{A}$ of Eq.~\ref{dense-subalg}, then we see that $\mathcal{A}$ is a $\C$-algebra and the elements $\{p_v, s_e, s_e^* : v \in E^0, e \in E^1\}$ satisfy the relations in Definition~\ref{Leavitt-def}.  Hence by the universal property of $L_K(E)$, there exists a homomorphism $\phi : L_\C(E) \to \mathcal{A}$ with $\phi(v) = p_v$, $\phi(e) = s_e$, and $\phi(e^*) = s_e^*$.  Furthermore, since $\{p_v, s_e, s_e^*: v \in E^0, e \in E^1\}$ generate $\mathcal{A}$ as an algebra, we see that $\phi$ is surjective.

We would like to show that $\phi$ is also injective, and hence an isomorphism.  An immediate idea of how to accomplish this is to use the universal property of $C^*(E)$ to obtain a homomorphism $\psi$ such that $\psi|_\mathcal{A} : \mathcal{A} \to L_\C(E)$ is an inverse for $\phi$.  However, there is a problem: in order to use the universal property of $C^*(E)$ to obtain the homomorphism $\psi$ we need to know that the generators $\{v, e, e^*:v \in E^0, e \in E^1\}$ of $L_K(E)$ sit inside a $C^*$-algebra.  Thus we need to show that $L_K(E)$ embeds into a $C^*$-algebra (or equivalently, that $L_K(E)$ embeds as a $*$-subalgebra of $\mathcal{B}(H)$, the bounded operators on a Hilbert space).  This is, in general, a difficult thing to show, and the author does not know of any elementary way to prove it for the algebra $\mathcal{A}$.

We can, however, show that the homomorphism $\phi$ is injective by using the Graded Uniqueness Theorem for Leavitt path algebras.  But, as one can see from \S \ref{GUT-sec}, this is a fairly nontrivial result.

\begin{theorem} \label{LC-embeds}
Let $E$ be a graph, and let $\phi : L_\C(E) \to C^*(E)$ be the canonical map onto the algebra $\mathcal{A}$ in Eq.~\ref{dense-subalg} obtained by the universal property of $L_\C(E)$.  Then $\phi$ is injective, and $L_\C(E)$ is isomorphic to a dense $*$-subalgebra of $C^*(E)$. 
\end{theorem}

\begin{proof}
By the universal property of $C^*(E)$ there exists a gauge action $\gamma : \T \to \operatorname{Aut} C^*(E)$ with $\gamma_z(p_v) = p_v$ for all $v \in E^0$ and $\gamma_z(s_e) = zs_e$ for all $e \in E^1$.  For $n \in \Z$ we may then define $\mathcal{A}_n := \{ a \in \mathcal{A} : \int_\T z^{-n} \gamma_z(a) \, dz = a \}$,  where the integration $dz$ is done with respect to normalized Haar measure on $\T$.  

We see that for an element $\lambda s_\alpha s_\beta^*$, we have $$\int_\T \lambda s_\alpha s_\beta^* \, dz = \begin{cases} \lambda s_\alpha s_\beta^* & \text{ if $|\alpha| - |\beta| = n$} \\ 0 & \text{ otherwise.} \end{cases}$$  Thus an element $x := \sum_{k=1}^N \lambda_k s_{\alpha_k} s_{\beta_k}^* \in \mathcal{A}$ is in $\mathcal{A}_n$ if an only if $|\alpha_k| - |\beta_k| = n$ for all $1 \leq k \leq N$.  One can then see that $\mathcal{A} = \bigoplus_{n \in \Z} \mathcal{A}_n$ as $\mathcal{A}$-modules.  Furthermore, if $x := \sum_{k=1}^M \lambda_k s_{\alpha_k} s_{\beta_k}^* \in \mathcal{A}_m$ and $y := \sum_{l=1}^N \kappa_l s_{\gamma_l} s_{\delta_l}^* \in \mathcal{A}_n$, we have that $$xy = \sum_{k,l} \eta_{k,l} s_{\mu_{k,l}} s_{\nu_{k,l}}^*,$$ where $|\mu_{k,l}| - |\nu_{k,l}| = |\alpha_k| - |\beta_k| + |\gamma_l| - |\delta_l| = m + n$.  Thus $xy \in \mathcal{A}_{m+n}$, and $\mathcal{A}$ is graded.  Since $\phi(v) = p_v \in \mathcal{A}_0$, $\phi(e) = s_e \in \mathcal{A}_1$, and $\phi(e^*) = s_e^* \in \mathcal{A}_{-1}$, we see that $\phi$ is a graded homomorphism.  Because we also have $\phi(v) = p_v \neq 0$ for all $v \in E^0$, it follows from Theorem~\ref{GUT} that $\phi$ is injective.
\end{proof}

This result has consequences for the construction of the graph $C^*$-algebra $C^*(E)$.  We recall this construction from \cite[Theorem~1.2]{KPR}:  

\begin{remark}[The Construction of $C^*(E)$] \label{graph-C*-construction}
Given a graph $E$, let $S_E := \{ (\alpha, \beta) : \alpha, \beta \in E^* \text{ and } r(\alpha) = r(\beta) \}$, and let $k_E$ be the space of complex-valued functions of finite support on $S_E$.  Then the set of point masses $\{\epsilon_\lambda : \lambda \in S_E \}$ forms a basis for the vector space $k_E$.  We define an associative multiplication on $S_E$ by setting
$$\epsilon_{(\alpha, \beta)} \epsilon_{(\gamma, \delta)} = \begin{cases} \epsilon_{(\alpha \gamma', \delta)} & \text{ if $\gamma = \beta \gamma'$} \\ \epsilon_{(\alpha , \delta)} & \text{ if $\beta = \gamma$} \\ \epsilon_{(\alpha, \beta'\delta)} & \text{ if $\beta = \gamma \beta'$} \\ 0 & \text{ otherwise.} \end{cases}$$  If we let $J$ be the ideal in $k_E$ generated by $\{\epsilon_{(v,v)}-\sum_{s(e)=v} \epsilon_{(e,e)} : v \in E^0_\textnormal{reg} \}$, then $k_E/J = L_\C(E)$. (Simply check that $k_E / J$ has the appropriate universal property.)

If we define 
\begin{align*}
\| a \|_0 := \sup \{ \| \pi(a) \| : &\text{ $\pi$ is a non-degenerate } \\ & \qquad \qquad \text{ $*$-representation of $k_E/J$ into $\mathcal{B} (H)$} \}
\end{align*}
then a standard argument shows that $\| \cdot \|_0$ is a well-defined, bounded semi-norm on $k_E /J$.  If we let $K := \{ a \in k_E/J : \| a \|_0 = 0 \}$, then $K$ is an ideal, and the completion $\overline{(k_E/J)/K}$ is a $C^*$-algebra.  One can then show that $\{ \epsilon_{(v,v)}, \epsilon_{(e, r(e))} : v \in E^0, e \in E^1 \}$ is a universal Cuntz-Krieger $E$-family generating $\overline{(k_E/J)/K}$, and thus $C^*(E) = \overline{(k_E/J)/K}$.
\end{remark}

Observe that the algebra $k_E/J$ in the above construction is the Leavitt path algebra $L_\C(E)$.  We see that $k_E$ is isomorphic to the path algebra generated by $E$ subject to relations (1)--(3) of Definition~\ref{Leavitt-def}, and $J$ is the ideal generated by the differences in relation (4) of Definition~\ref{Leavitt-def}, which is precisely how $L_\C(E)$ is constructed.  Therefore, in order for the canonical map from $L_\C$ onto the subalgebra $\mathcal{A} = \{ a + K : a \in k_E/J \}$ to be injective, we would need $K$ to be zero.

If one looks at $\| \cdot \|_0$, it is easy to see that $\| \cdot \|_0$ is a semi-norm.  However, it is not clear at all that $\| \cdot \|_0$ is a norm.  (To show this one would need to show that for every $a \in k_E/J$ there is a $*$-representation $\pi : k_E / J \to \mathcal{B}(H)$ with $\pi(a) \neq 0$.)  Yet, because we know from Theorem~\ref{LC-embeds} that the map $\phi : L_\C(E) \to C^*(E)$ embeds $L_\C(E)$ onto $\mathcal{A}$, it must be the case that $K = \{ 0 \}$ and $\| \cdot \|_0$ is a norm.  This is a fact that is likely to be surprising to most $C^*$-algebraists!  We summarize these consequences in the following two corollaries.

\begin{corollary} \label{graph-C*-construct-cor}
The ideal $K$ in the construction of the graph $C^*$-algebra described in Remark~\ref{graph-C*-construction} is equal to $\{ 0 \}$, and $C^*(E) = \overline{k_E/J}$.
\end{corollary}

\begin{corollary} \label{norm-LC-in-C*}
Let $E$ be a graph and let $L_\C(E)$ be the associated Leavitt path algebra with coefficients in $\C$.  If we define $\| \cdot \|_0 : L_\C(E) \to \C$ by 
\begin{align*}
\| a \|_0 := \sup \{ \| \pi(a) \| : &\text{ $\pi$ is a non-degenerate } \\ & \qquad \qquad \text{ $*$-representation of $L_\C(E)$ into $\mathcal{B} (H)$} \}
\end{align*}
then $\| \cdot \|_0$ is a norm on $L_\C(E)$ and $C^*(E) = \overline{L_\C(E)}^0$, where $\overline{L_\C(E)}^0$ denotes the completion of $L_\C(E)$ with respect to the norm $\| \cdot \|_0$.
\end{corollary}

Since $L_\C(E)$ may be viewed as a $*$-subalgebra of $C^*(E)$, we can take an ideal in $L_\C(E)$ and form its closure to obtain an ideal in $C^*(E)$.  As shown in the following proposition, this gives a correspondence between graded ideals of $L_\C(E)$ and gauge-invariant ideals of $C^*(E)$.

\begin{proposition}
Let $E$ be a graph.  Identify $L_\C(E)$ as a $*$-subalgebra of $C^*(E)$ as described in Corollary~\ref{norm-LC-in-C*}, and for an ideal $I$ of $L_\C(E)$ let $\overline{I}$ denote the closure of $I$ in $C^*(E)$.  Then the map $$I \mapsto \overline{I}$$ is a lattice isomorphism from the lattice of graded ideals of $L_\C(E)$ onto the lattice of of gauge-invariant ideals of $C^*(E)$, with inverse given by $J \mapsto J \cap L_\C(E)$.

Moreover, when $E$ satisfies Condition~(K) all ideals of $L_\C(E)$ are graded, all ideals of $C^*(E)$ are gauge-invariant, and the map $I \mapsto \overline{I}$ is a lattice isomorphism from the lattice of ideals of $L_\C(E)$ onto the lattice of ideals of $C^*(E)$.
\end{proposition}

\begin{proof}
For an admissible pair $(E,S)$ of the graph $E$, let $I_{(H,S)}$ denote the ideal in $L_\C(E)$ generated by $\{ v : v \in H\} \cup \{ v^H : v \in S \}$.  In addition, identify the elements $v$ and $v^H := v - \sum_{s(e)} ee^*$ in $L_\C(E)$ with the elements $p_v$ and $p_v^H := p_v -\sum_{s(e)=v} s_es_e^*$ in $C^*(E)$, and let $J_{(H,S)}$ be the ideal in $C^*(E)$ generated by $\{ p_v : v \in H\} \cup \{ p_v^H : v \in S \}$.

We shall prove that $\overline{I_{(H,S)}} = J_{(H,S)}$.  To begin, we note that $\overline{I_{(H,S)}}$ is an ideal in $C^*(E)$ containing $\{ p_v : v \in H\} \cup \{ p_v^H : v \in S \}$.  Therefore $J_{(H,S)} \subseteq \overline{I_{(H,S)}}$.  Furthermore, if $x \in \overline{I_{(H,S)}}$, then by Lemma~\ref{lem-ideal-span} $x = \lim_{n} \sum_{k=1}^{N_n} \lambda_{k,n} \alpha_{k,n} \beta_{k,n}^*$ for $\alpha_{n,k}, \beta_{n,k} \in E^*$.  Using the identification mentioned above, we have that $x = \lim_{n} \sum_{k=1}^{N_n} \lambda_{k,n} s_{\alpha_{k,n}} s_{\beta_{k,n}}^*$.  It follows from the first paragraph of \cite[p.6]{BHRS} that $x \in J_{(H,S)}$.  Therefore, $\overline{I_{(H,S)}} = J_{(H,S)}$.

It follows from Theorem~\ref{ideal-structure-theorem}(\ref{one}) that $(H,S) \mapsto I_{(H,S)}$ is a lattice isomorphism from the admissible pairs of $E$ onto the graded ideals of $L_\C(E)$, and it follows from \cite[Theorem~3.6]{BHRS} that $(H,S) \mapsto J_{(H,S)}$ is a lattice isomorphism from the admissible pairs of $E$ onto the gauge-invariant ideals of $C^*(E)$.  If we compose the inverse of the first map with the second map, we see that $I_{(H,S)} \mapsto \overline{I_{(H,S)}} = J_{(H,S)}$ is a lattice isomorphism from the graded ideals of $L_\C(E)$ onto the gauge-invariant ideals of $C^*(E)$.  Furthermore, it is straightforward to show that the inverse is given by $J \mapsto J \cap L_\C(E)$.

Moreover, when $E$ satisfies Condition~(K), it follows from Theorem~\ref{K-iff-ideals-graded} that every ideal of $L_\C(E)$ is graded, and it follows from \cite[Corollary~3.8]{BHRS} (or \cite[Theorem~3.5]{DT1}) that every ideal of $C^*(E)$ is gauge invariant.  Hence when $E$ satisfies Condition~(K), then map $I_{(H,S)} \mapsto \overline{I_{(H,S)}}$ is a map from the lattice of all ideals of $L_\C(E)$ onto the lattice of all ideals of $C^*(E)$.
\end{proof}

\begin{remark}
We mention that the above result does not hold if we extend the map $I \mapsto \overline{I}$ to ideals that are not graded.  If $E$ is the graph

$ $

$$
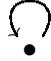
\xymatrix{
\bullet \ar@(ur,ul) \\
}
$$
consisting of a single vertex and a single edge, then $L_\C(E) \cong \C [ x, x^{-1}]$ and $C^*(E) \cong C(\T)$.  The inclusion of $L_\C(E) \hookrightarrow C^*(E)$ identifies $L_\C(E)$ with the finite Laurent polynomials $p(z, z^{-1})$ on $\T$.

If we view $\T = [0, 2 \pi)$ and let $[a,b]$, with $a \neq b$, be an interval of $\T$, we can let  $J := \{ f \in C(\T) : f|_{[a,b]} = 0 \}$.  Then $J$ is an ideal of $C^*(E)= C(\T)$, but $J \cap L_\C(E) = \{0 \}$ since any nonzero polynomial can only vanish at a finite number of points.  Hence the map $I \mapsto \overline{I}$ is not surjective when extended to all ideals.

Furthermore, if we let $a = 1 + z + z^3$, then as shown in Remark~\ref{non-selfadjoint-ideal-rem}, the ideal $I_1 := \langle a \rangle$ in $L_\C(E)$ is not self-adjoint and in particular $a^* = 1 + z^{-1} + z^{-3} \notin I_1$.  However, $\overline{I_1}$ is an ideal in the $C^*$-algebra $C^*(E)$, and thus $\overline{I_1}$ is self adjoint.  Therefore if we let $I_2 := \langle a^* \rangle$, then $\overline{I_1} = \overline{I_2}$ even though $I_1 \neq I_2$.  Hence the map $I \mapsto \overline{I}$ is not injective when extended to all ideals.
\end{remark}

\noindent \textbf{Acknowledgements:} The author thanks Enrique Pardo for many useful suggestions and comments regarding a preliminary draft of this article.  The author also thanks Pere Ara for his comments and for pointing out a mistake in a prior version of the proof of Lemma~\ref{F-G-comm-diagram}.

\end{document}